\documentclass{tglat2e}
\usepackage{graphics}
\usepackage{amssymb}
\usepackage{latexsym,bm}
\usepackage{amsmath,amsfonts,amssymb,mathrsfs}
\usepackage[all]{xy}
\usepackage{amscd}

\DeclareMathOperator{\End}{End} \DeclareMathOperator{\Hom}{Hom}
\DeclareMathOperator{\Bid}{Id} \DeclareMathOperator{\newid}{id}
\DeclareMathOperator{\chf}{ch} \DeclareMathOperator{\BBd}{Bd}
\DeclareMathOperator{\Ker}{Ker} 
\DeclareMathOperator{\rank}{rank}

\DeclareMathOperator{\Ext}{Ext}
\DeclareMathOperator{\Ind}{Ind}
\DeclareMathOperator{\Image}{Im} \DeclareMathOperator{\soc}{soc}
\newcommand{\mHH}{\mathcal H}
\DeclareMathOperator{\sign}{sign}

\newcommand{\mbb}{\mathfrak{B}}
\newcommand{\mmZ}{\mathbb Z}
\newcommand{\mmA}{\mathscr A}

\newcommand{\mmU}{\mathbb U}
\newcommand{\mmC}{\mathbb C}
\newcommand{\mmQ}{\mathbb Q}

\def\mui{\underline{i}}
\def\muj{\underline{j}}

\newcommand{\lam}{\lambda}
\newcommand{\eps}{\varepsilon}

\newcommand{\bBS}{\mathfrak S}
\newcommand{\mft}{\mathfrak t}

\theoremstyle{plain}
\newtheorem{theorem}{Theorem}
\newtheorem{lemma}[theorem]{Lemma}
\newtheorem{corollary}[theorem]{Corollary}
\newtheorem{proposition}[theorem]{Proposition}

\theoremstyle{definition}
\newtheorem{definition}[theorem]{Definition}

\theoremstyle{remark}
\newtheorem{remark}[theorem]{Remark}
\numberwithin{theorem}{section}

\begin{document}

\title[Dual partially harmonic tensors]
{Dual partially harmonic tensors and Brauer--Schur--Weyl duality}

\authors{J.~Hu\thanks{Supported by an ARC discovery grant, an NSF of China (No. 10771014), the Scientific Research
Foundation for the Returned Overseas
Chinese Scholars by the State Education Ministry and the Basic Research Foundation of BIT.}
\address Department of Mathematics\\
Beijing Institute of Technology\\
Beijing, 100081, P.R. China\\
\quad\,\, \&\\
School of Mathematics and Statistics\\
University of Sydney\\
NSW 2006, Australia \email junhu303@yahoo.com.cn}
\received{???}
\accepted{???}

\maketitle
\begin{abstract} Let $V$ be a $2m$-dimensional symplectic vector space over an algebraically closed field $K$.
Let $\mbb_n^{(f)}$ be the two-sided ideal of the  Brauer algebra
$\mbb_n(-2m)$ over $K$ generated by $e_1e_3\cdots e_{2f-1}$, where
$0\leq f\leq [n/2]$. Let $\mathcal{HT}_{f}^{\otimes n}$ be the
subspace of partially harmonic tensors of valence $f$ in $V^{\otimes
n}$. In this paper, we prove that $\dim\mathcal{HT}_f^{\otimes n}$
and $\dim\End_{KSp(V)}\Bigl(V^{\otimes n}/V^{\otimes
n}\mbb_n^{(f)}\Bigr)$ are both independent of $K$, and the natural
homomorphism from $\mbb_n(-2m)/\mbb_n^{(f)}$ to
$\End_{KSp(V)}\Bigl(V^{\otimes n}/V^{\otimes n}\mbb_n^{(f)}\Bigr)$
is always surjective. We show that $\mathcal{HT}_{f}^{\otimes n}$
has a Weyl filtration and is isomorphic to the dual of $V^{\otimes
n}\mbb_n^{(f)}/V^{\otimes n}\mbb_n^{(f+1)}$ as a
$Sp(V)$-$(\mbb_n(-2m)/\mbb_n^{(f+1)})$-bimodule. We obtain a
$Sp(V)$-$\mbb_n$-bimodules filtration of $V^{\otimes n}$ such that
each successive quotient is isomorphic to some $\nabla(\lam)\otimes
z_{g,\lam}\mbb_n$ with $\lam\vdash n-2g$, $\ell(\lam)\leq m$ and
$0\leq g\leq [n/2]$, where $\nabla(\lam)$ is the co-Weyl module
associated to $\lam$ and $z_{g,\lam}$ is an explicitly constructed
maximal vector of weight $\lam$. As a byproduct, we show that each
right $\mbb_n$-module $z_{g,\lam}\mbb_n$ is integrally defined and
stable under base change.
\end{abstract}

\section{Introduction}

Let $m, n\in\mathbb{N}$. Let $K$ be an algebraically closed field
and $V$ a $2m$-dimensional symplectic vector space over $K$. The
symplectic group $Sp(V)$ acts naturally on $V$ from the left hand
side, and hence on the $n$-tensor space $V^{\otimes n}$. Let
$\mbb_n:=\mbb_n(-2m)$ be the Brauer algebra over $K$ with
generators $s_1,\cdots,s_{n-1},e_1,$ $\cdots,e_{n-1}$ and parameter
$-2m\cdot 1_{K}$ (see \ref{dfn1} for their definitions). There is a right action of $\mbb_n$ on $V^{\otimes
n}$ which commutes with the left action of $Sp(V)$. Let
$\varphi_K,\psi_K$ be the following natural $K$-algebra
homomorphisms:
$$
\varphi_K:\,(\mbb_n)^{\rm op}\rightarrow\End_{KSp(V)}\bigl(V^{\otimes
n}\bigr),\quad\, \psi_K:\,
KSp(V)\rightarrow\End_{\mbb_n}\bigl(V^{\otimes n}\bigr).
$$
For any positive integer $k$, a composition of $k$ is a sequence of
non-negative integers $\lam=(\lam_1,\lam_2,\cdots)$ with $\sum_{i\geq 1}\lam_i=k$.  We use $\ell(\lam)$ to denote the largest integer $t$ such that $\lam_t\neq 0$. A composition $\lam=(\lam_1,\lam_2,\cdots)$ of $k$ is
said to be a partition if $\lam_1\geq\lam_2\geq\cdots$; in that
case, we write $\lam\vdash k$. The following results
are referred as Brauer--Schur--Weyl duality.

\begin{theorem} \label{thm11} {\rm (\cite{B}, \cite{DP}, \cite{DDH})}
\label{ddh1} 1) The natural left
action of $Sp(V)$ on $V^{\otimes n}$ commutes with the right action
of $\mbb_n$;

2) both $\varphi_K$ and $\psi_K$ are
surjective;

3) if $m\geq n$ then $\varphi_K$ is injective, and hence an
isomorphism;

4) if $K=\mathbb{C}$, then there is a decomposition of $V^{\otimes n}$ as a direct sum of
irreducible $\mmC Sp(V)$-$\mbb_n$-bimodules:
$$
V^{\otimes n}=\bigoplus_{f=0}^{[n/2]}\bigoplus_{\substack{\lam\vdash n-2f\\
\ell(\lam)\leq m}}\Delta({\lam})\otimes D^{(f,\lam)},$$ where
$\Delta({\lam})$ (respectively, $D^{(f,\lam)}$) denotes the
irreducible $Sp(V)$-module (respectively, the irreducible
$\mbb_n$-module) corresponding to $\lam$ (respectively, corresponding
to $(f,\lam)$).
\end{theorem}

There is a variant of the above Brauer--Schur--Weyl duality as we
shall describe. Let $\mbb_n^{(1)}$ be the two-sided ideal of $\mbb_n$
generated by $e_1$. We set $$ \quad W_{1,n}:=\bigl\{v\in V^{\otimes
n}\bigm|vx=0,\forall\,x\in \mbb_n^{(1)}\bigr\}.
$$

\begin{definition} \label{w1n} We call $W_{1,n}$ the subspace of harmonic tensors or traceless tensors.
\end{definition}

Note that our definition of harmonic tensors looks slightly different
from that given in \cite{DS} and \cite[\S 10.2.1]{GW}. The two
definitions are reconciled in Corollary \ref{2dfn}. \smallskip

If $K=\mathbb{C}$, then we shall write $V_{\mmC}, W^{\mmC}_{1,n}$
instead of $V, W_{1,n}$ in order to emphasize the base field. We have
that $\mbb_n/\mbb_n^{(1)}\cong K\bBS_n$. The right action of $\mbb_n$ on
$V^{\otimes n}$ gives rise to a right action of $K\bBS_n$ on $W_{1,n}$.

\begin{theorem} {\rm (\cite[(10.2.7)]{GW}, \cite{W})} \label{GW1} The natural left action of $Sp(V)$
on $W_{1,n}$ commutes with the right action of $K\bBS_n$. If $K=\mmC$,
then there is a decomposition of $W^{\mmC}_{1,n}$ as a direct sum of
irreducible $\mmC Sp(V)$-$\mmC\bBS_n$-bimodules:
$$
W^{\mmC}_{1,n}=\bigoplus_{\substack{\lam\vdash n\\
\ell(\lam)\leq m}}\Delta_{\mmC}({\lam})\otimes S_{\mmC}^{\lam},$$ where
$\Delta_{\mmC}({\lam})$ (respectively, $S_{\mmC}^{\lam}$) denotes the
irreducible $\mmC Sp(V)$-module (respectively, the irreducible
$\mmC\bBS_n$-module) corresponding to $\lam$.
\end{theorem}

As before, we have two natural $K$-algebra homomorphisms:
$$
\varphi^{(1)}_{K}:\,(K\bBS_n)^{\rm op}\rightarrow\End_{KSp(V)}\bigl(W_{1,n}\bigr),\quad\,
\psi^{(1)}_{K}:\, KSp(V)\rightarrow\End_{K\bBS_n}\bigl(W_{1,n}\bigr).
$$
In \cite{DS}, De Concini and Strickland proved that $\dim W_{1,n}$
is independent of the field $K$ and $\varphi^{(1)}_K$ is always surjective.
Furthermore, they showed that if $m\geq n$, then $\varphi^{(1)}_K$
is an isomorphism. Their proof makes use of the previous results in
\cite{De} and \cite{DP} on multilinear invariants of a variety and
symplectic standard tableaux which eventually relies on some
algebro-geometric arguments. In \cite{Ma2}, using the theory of
rational representations of symplectic group, Maliakas proved that
$W_{1,n}^{\ast}$ has a good filtration whenever $m\geq n$ and he
claimed that it is true for arbitrary $m$.\smallskip

The starting point of this paper is, on one hand, to generalize
the above duality to the case of partially harmonic tensors of
arbitrary valence $f$, and on the other hand, to provide a
self-contained and purely representation-theorietic approach which
makes it possible to work also in the quantized case\footnote{At the
moment, there are still a few obstacles (e.g., Lemma \ref{keylem1})
which prevent us from generalizing the main results of this paper
to the quantized case. The main difficulty lies in that there are various choices of tangle which can specialize to the same Brauer diagram and it
also becomes much harder to describe the action of a tangle on tensor space in a simple combinatorial way (as in the Brauer algebra case).}. We are mostly interested in the
non-semisimple case. To describe our main results, we need some more
notations and definitions. For each integer $0\leq f\leq [n/2]$, let
$\mbb_n^{(f)}$ be the two-sided ideal of $\mbb_n$ generated by
$e_1e_3\cdots e_{2f-1}$. By convention, $\mbb_n^{(0)}=\mbb_n$ and $\mbb_n^{([n/2]+1)}=\{0\}$. Set
$$
\mathcal{HT}_f^{\otimes n}:=\bigl\{v\in V^{\otimes
n}\mbb_n^{(f)}\bigm|vx=0, \forall\, x\in\mbb_n^{(f+1)}\bigr\}.
$$
Following \cite[(10.3.1)]{GW}, we call $\mathcal{HT}_f^{\otimes n}$
the space of {\it partially harmonic tensors of valence $f$}. Classically, these spaces play an important role in the study of the $\mathbb{C}Sp(V_{\mmC})$-module structure on the $n$-tensor spaces $V_{\mathbb{C}}^{\otimes n}$ (cf. \cite{GW}). One of our motivation for studying them is to try to use them to construct some natural (integral) quotient of the symplectic Schur algebras. By
\cite[(10.3.14)]{GW}, it is easy to see that if $K=\mathbb{C}$ then $\mathcal{HT}_f^{\otimes n}\cong V_{\mmC}^{\otimes
n}\mbb_n^{(f)}/V_{\mmC}^{\otimes n}\mbb_n^{(f+1)}$
as a $\mmC$-linear space. In particular,
$W_{1,n}^{\mmC}=\mathcal{HT}_0^{\otimes n}\cong V_{\mmC}^{\otimes
n}/V_{\mmC}^{\otimes n}\mbb_n^{(1)}$ as a $\mmC$-linear space. \smallskip

Note that the space $V$ can be defined over arbitrary field (and even over $\mmZ$). So we can consider the spaces $\mathcal{HT}_f^{\otimes
n}$,  $V^{\otimes n}\mbb_n^{(f)}/V^{\otimes n}\mbb_n^{(f+1)}$ and
$V^{\otimes n}/V^{\otimes n}\mbb_n^{(f)}$ over an arbitrary
field $K$. All the results we obtain in this paper are valid over these more general ground fields.
However, we shall always assume that $K$ is algebraically closed whenever notions and results from algebraic groups theory (e.g., good filtration, Weyl filtration) are needed. The transition between an arbitrary field and its algebraic closure usually follows from some standard arguments in commutative algebras. Note that $V^{\otimes n}/V^{\otimes n}\mbb_n^{(f)}$ is a
$Sp(V)$-$\bigl(\mbb_n/\mbb_n^{(f)}\bigr)$-bimodule. We use
$\varphi_{f,K}$ to denote the following natural homomorphism: $$
\varphi_{f,K}:\,\mbb_n/\mbb_n^{(f)}\rightarrow
\End_{KSp(V)}\Bigl(V^{\otimes n}/V^{\otimes n}\mbb_n^{(f)}\Bigr).$$

Let $\delta_{ij}$ denote the value of the usual Kronecker delta. For
each integer $1\leq i\leq 2m$, we set $i':=2m+1-i$. We fix an
ordered basis $\bigl\{v_1,v_2,\cdots,v_{2m}\bigr\}$ of $V$ such that
$$ \langle v_i, v_{j}\rangle=0=\langle v_{i'},
v_{j'}\rangle,\,\,\,\langle v_i, v_{j'}\rangle=\delta_{ij}=-\langle
v_{j'}, v_{i}\rangle,\quad\forall\,\,1\leq i, j\leq m. $$ We use
$V_{\mmZ}$ to denote the free $\mmZ$-submodule of $V_{\mmC}$ generated by
$v_1, \cdots, v_{2m}$. For any commutative $\mmZ$-algebra $R$, we set
$V_{R}:=V_{\mmZ}\otimes_{\mmZ}R$. The Brauer algebra $\mbb_n$ can also be
defined over $R$ and we denote it
by $\mbb_n^{R}$. To simplify notations, the two-sided ideal of
$\mbb_n^{R}$ generated by $e_1e_3\cdots e_{2f-1}$ will be still
denoted by $\mbb_n^{(f)}$. The main results in this paper are the
following theorems and corollaries.

\begin{theorem} \label{mainthm0} For each integer $1\leq f\leq [n/2]$,

1) $V_{\mmZ}^{\otimes n}\mbb_n^{(f)}$ is a pure $\mmZ$-submodule of
$V_{\mmZ}^{\otimes n}$, equivalently, $V_{\mmZ}^{\otimes
n}/V_{\mmZ}^{\otimes n}\mbb_n^{(f)}$ is a free $\mmZ$-module;

2) both $V^{\otimes n}\mbb_n^{(f)}$ and $V^{\otimes n}/V^{\otimes n}\mbb_n^{(f)}$ are stable under base change, i.e., for
any commutative $\mmZ$-algebra $R$, the canonical maps
$$
V_{\mmZ}^{\otimes n}\mbb_n^{(f)}\otimes_{\mmZ}R\rightarrow V_{R}^{\otimes
n}\mbb_n^{(f)},\,\,\,V_{\mmZ}^{\otimes n}/V_{\mmZ}^{\otimes
n}\mbb_n^{(f)}\otimes_{\mmZ}R\rightarrow V_{R}^{\otimes
n}/V_{R}^{\otimes n}\mbb_n^{(f)}$$ are isomorphisms. In particular,
the character formulae of the left $Sp(V)$-modules $V^{\otimes
n}\mbb_n^{(f)}, V^{\otimes n}/V^{\otimes n}\mbb_n^{(f)}$ are both
independent of the field $K$.
\end{theorem}

\begin{theorem}\label{mainthm1} Let $K$ be an algebraically closed field. For each integer $f$ with $1\leq f\leq [n/2]$, both $V^{\otimes n}/V^{\otimes n}\mbb_n^{(f)}$ and $V^{\otimes n}\mbb_n^{(f)}$ have a good filtration as $Sp(V)$-modules.
\end{theorem}

\begin{corollary}\label{maincor1} Let $K$ be an algebraically closed field. Let $0\leq f\leq [n/2]$. Then the dimension of
$V^{\otimes n}\mbb_n^{(f)}/V^{\otimes n}\mbb_n^{(f+1)}$ is independent of $K$, and there is a
$Sp(V)$-$(\mbb_n/\mbb_n^{(f+1)})$-bimodule isomorphism:
$$V^{\otimes n}\mbb_n^{(f)}/V^{\otimes n}\mbb_n^{(f+1)}\cong \Bigl(\mathcal{HT}_f^{\otimes n}\Bigr)^{\ast}. $$
\end{corollary}

In particular, the dimension of $\mathcal{HT}_f^{\otimes n}$ is
independent of $K$ too.

\begin{corollary} \label{maincor15} Let $K$ be an algebraically closed field. For each $0\leq f\leq [n/2]$, the $Sp(V)$-module
$V^{\otimes n}\mbb_n^{(f)}/V^{\otimes n}\mbb_n^{(f+1)}$
always has a good filtration and the $Sp(V)$-module
 $\mathcal{HT}_f^{\otimes n}$ always has a Weyl filtration.
\end{corollary}

\begin{theorem}\label{mainthm2} Let $K$ be an algebraically closed field. For each integer $f$ with $1\leq f\leq [n/2]$,

1) the dimension of the endomorphism algebra
$\End_{KSp(V)}\Bigl(V^{\otimes n}/V^{\otimes n}\mbb_n^{(f)}\Bigr)$ is
independent of $K$;

2) $\varphi_{f,K}(\mbb_n/\mbb_n^{(f)})=\End_{KSp(V)}\Bigl(V^{\otimes
n}/V^{\otimes n}\mbb_n^{(f)}\Bigr)$.
\end{theorem}

The proof of the above results are {\it completely self-contained} and use purely representation-theorietic knowledge.  As a consequence of these theorems and corollaries, we recover
and extend the previously mentioned results of De Concini and
Strickland \cite{DS} and the result of Maliakas \cite{Ma2}. In the way of our proof, we also obtain the following result, which seems of independent interest.

\begin{theorem}\label{mainthm3} 1) As a $Sp(V)$-$\mbb_n$-bimodule, $V^{\otimes n}$ has a filtration such that each successive quotient is isomorphic to some $\nabla(\lam)\otimes z_{g,\lam}\mbb_n$ with $\lam\vdash n-2g$, $\ell(\lam)\leq m$ and
$0\leq g\leq [n/2]$, where
$\nabla(\lam)$ is the co-Weyl module associated to $\lam$ and
$z_{g,\lam}$ is a maximal vector of weight $\lam$ (see \ref{Zdfn}
for its definition);

2) for any partition $\lam $ of $n-2g$ with $0\leq g\leq [n/2]$ and $\ell(\lam)\leq m$ and any commutative
$\mmZ$-algebra $R$, the canonical map $$
z_{g,\lam}\mbb_n^{\mmZ}\otimes_{\mmZ}R\rightarrow z_{g,\lam}\mbb_n^{R}
$$
is always an isomorphism. In particular, the dimension of $z_{g,\lam}\mbb_n$ is independent of $K$.
\end{theorem}
In fact, if $K=\mathbb{C}$ then $z_{g,\lam}\mbb_n^{\mmC}$ is a simple
right $\mbb_n^{\mmC}$-module. Therefore, for any field $K$, the
dimension of $z_{g,\lam}\mbb_n$ is always equal to the number of
$(-2m)$-permissible up-down tableaux of shape $\lam'$ and length $n$
(cf. \cite[Theorem 1.1, Theorem 1.2]{HY} and \cite{Wz}), where $\lam'$ denotes the conjugate of $\lam$.\smallskip

The paper is organized as follows. In Section 2 we recall some
basic knowledge about Brauer algebras and their actions on
$n$-tensor spaces. In particular, we show that our Definition
\ref{w1n} of harmonic tensors coincides with that given in \cite{DS}
and \cite[\S 10.2.1]{GW}. In Section 3, we give the proof of Theorem
\ref{mainthm0}, \ref{mainthm1} and Corollary \ref{maincor1}. The main idea
is to show that $V^{\otimes n}\mbb_n^{(f)}$ can be identified
with the image of $V^{\otimes n}$ under a truncation functor
$\mathcal{O}_{\pi_f}$ associated with a saturated set $\pi_f$ of
dominant weights. The proof makes use of the main result obtained in
\cite{DDH}, some results on Weyl filtration (resp., good filtration)
and a key Lemma \ref{Hu2lem}. As a consequence, we prove the first part of
Theorem \ref{mainthm3}. Section 4 is devoted to the proof of Lemma
\ref{Hu2lem}. The proof relies on Lusztig's theory of canonical bases and based modules. As a result, we get the
second part of Theorem \ref{mainthm3}. In Section 5 we prove Theorem \ref{mainthm2}, which
gives one side of the Brauer-Schur-Weyl duality between
$\mbb_n/\mbb_n^{(f)}$ and $KSp(V)$ on $V^{\otimes n}/V^{\otimes
n}\mbb_n^{(f)}$. We conjecture that the other side of this duality is also true. \smallskip

\section{Preliminaries}

Let $m,n\in\mathbb{N}$. The Brauer algebra $\mbb_n$ with parameter
$-2m$ and size $n$ was first introduced by Richard Brauer (see
\cite{B}) when he studied how the $n$-tensor space
$V_{\mmC}^{\otimes{n}}$ decomposes into irreducible modules over
$Sp(V_{\mmC})$. In his language, $\mbb_n$ was defined as the
$\mmC$-linear space with a basis being the set $\BBd_n$ of all the
Brauer $n$-diagrams. By definition, a Brauer $n$-diagram is a
diagram with specific $2n$ vertices arranged in two rows of $n$
each, the top and  bottom rows, and exactly $n$ edges such that
every vertex is joined to another vertex (distinct from itself) by
exactly one edge. We label the vertices in each row of a Brauer
$n$-diagram by the integers $1,2,\cdots,n$ from left to right. The
multiplication of two Brauer $n$-diagrams is defined as follows. We
compose two diagrams $D_1, D_2$ by identifying the bottom row of
$D_1$ with the top row of $D_2$ (such that the vertex $i$ in the
bottom row of $D_1$ is identified with the vertex $i$ in the top row
of $D_2$ for each $1\leq i\leq n$). The result is a graph, with a
certain number, $n(D_1,D_2)$, of interior loops. After removing the
interior loops and the identified vertices, retaining the edges and
remaining vertices, we obtain a new Brauer $n$-diagram $D_1\circ
D_2$, the composite diagram. Then we define $D_1\cdot
D_2=(-2m)^{n(D_1,D_2)}D_1\circ D_2$. In general, the multiplication
of two elements in $\mbb_n$ is given by the linear extension of a
product defined on diagrams. For example, let $D$ be the following
Brauer $5$-diagram.
\medskip
\begin{center}
\includegraphics{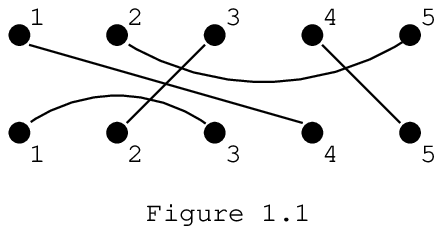}
\end{center}
\medskip
Let $D'$ be the following Brauer $5$-diagram.
\medskip
\begin{center}
\includegraphics{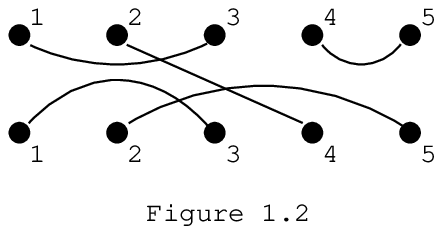}
\end{center}
\medskip
Then $DD'$ is equal to
\medskip
\begin{center}
\scalebox{0.55}[0.55]{\includegraphics{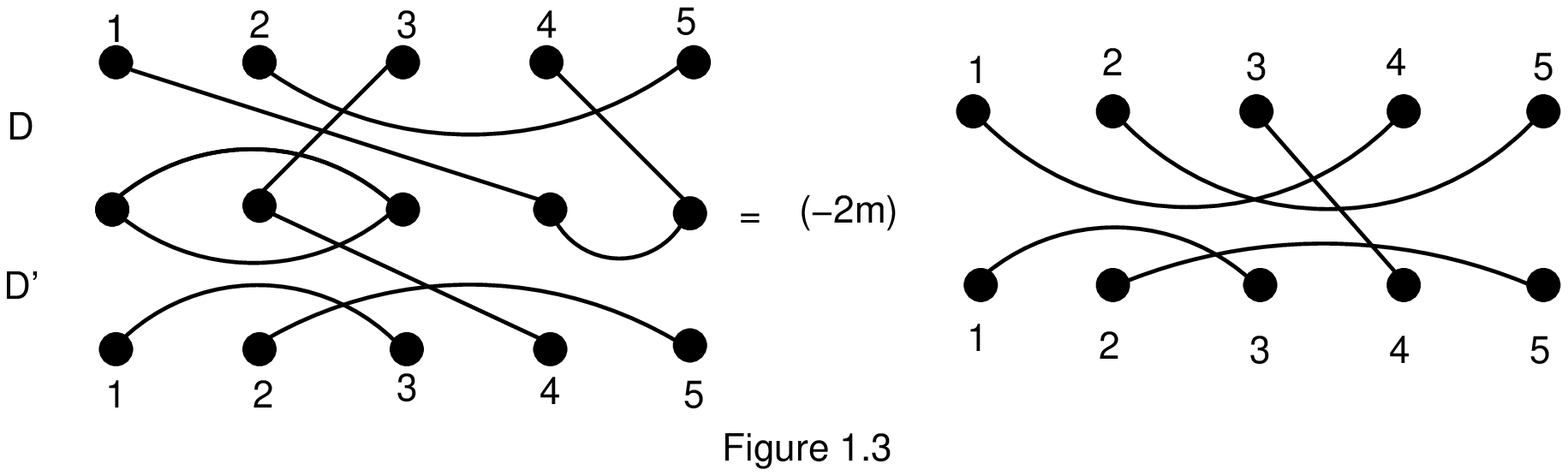}}
\end{center}
\medskip

Alternatively, one can define the Brauer algebra using generators
and relations.

\begin{definition} \label{dfn1} Let $K$ be a field. The Brauer algebra $\mbb_n$ over $K$ is a
unital associative $K$-algebra with canonical generators
$s_1,\cdots,s_{n-1},e_1,\cdots,e_{n-1}$ and relations (see
\cite{E}):
$$\begin{matrix}s_i^2=1,\,\,e_i^2=-2m e_i,\,\,e_is_i=e_i=s_ie_i,
\quad\forall\,1\leq i\leq n-1,\\
s_is_j=s_js_i,\,\,s_ie_j=e_js_i,\,\,e_ie_j=e_je_i,\quad\forall\,1\leq
i<j-1\leq n-2,\\ s_is_{i+1}s_i=s_{i+1}s_is_{i+1},\,\,
e_ie_{i+1}e_i=e_i,\,\, e_{i+1}e_ie_{i+1}=e_{i+1},\,\,\forall\,1\leq
i\leq n-2,\\
s_ie_{i+1}e_i=s_{i+1}e_i,\,\,e_{i+1}e_is_{i+1}=e_{i+1}s_i,\quad\forall\,1\leq
i\leq n-2.\end{matrix}
$$
\end{definition}

Replacing $K$ by any commutative $\mmZ$-algebra $R$, we can define the
Brauer algebra $\mbb_n^{R}$ over $R$ in a similar way. The algebra $\mbb_n^{R}$ is a free $R$-module with rank
$(2n-1)\cdot (2n-3)\cdots 3\cdot 1$, and the canonical map
$\mbb_n^{\mmZ}\otimes_{\mmZ}R\rightarrow\mbb_n^{R}$ is an isomorphism for
any commutative $\mmZ$-algebra $R$.\smallskip

The two definitions of the Brauer algebra $\mbb_n$ can be identified
as follows. Let $i$ be an integer with $1\leq i\leq n-1$. The
generator $s_i$ corresponds to the Brauer $n$-diagram with edges
connecting vertices $i$ (respectively, $i+1$) on the top row with
$i+1$ (respectively, $i$) on the bottom row, and all other edges are
vertical, connecting vertices $k$ on the top and bottom rows for all
$k\neq i,i+1$. The generator $e_i$ corresponds to the Brauer
$n$-diagram with horizontal edges connecting vertices $i,i+1$ on the
top and bottom rows, and all other edges are vertical, connecting
vertices $k$ on the top and bottom rows for all $k\neq i,i+1$. Note that the subalgebra of $\mbb_n^{\mmZ}$
generated by $s_1,s_2,\cdots,s_{n-1}$ is isomorphic to the group
algebra of the symmetric group $\bBS_{n}$ over $\mmZ$.\smallskip

The Brauer algebra was studied in a number of literatures, see
\cite{DDH}, \cite{DH}, \cite{E}, \cite{GL}, \cite{HW1},
\cite{HW2}, \cite{Hu1}. We are mainly interested in their actions on $n$-tensor
space $V^{\otimes n}$ and related Schur--Weyl dualities involving
symplectic groups. From now on let $K$ be an algebraically closed field and
$V$ a $2m$-dimensional $K$-vector space equipped with a
non-degenerate skew-symmetric bilinear form $\langle\,,\rangle$. Let
$Sp(V)$ be the corresponding symplectic group.  Recall that for
each integer $1\leq i\leq 2m$, $i':=2m+1-i$. Let $\{v_i\}_{i=1}^{2m}$ be a $K$-basis of $V$
such that $$\langle v_i, v_{j}\rangle=0=\langle v_{i'},
v_{j'}\rangle,\,\,\,\langle v_i, v_{j'}\rangle=\delta_{ij}=-\langle
v_{j'}, v_{i}\rangle,\quad\forall\,\,1\leq i, j\leq m. $$
For each integer $1\leq i\leq 2m$, we define $$ v_i^{\ast}=\begin{cases}
v_{i'}, &\text{if $1\leq i\leq m$;}\\
-v_{i'}, &\text{if $m+1\leq i\leq 2m$.}
\end{cases}
$$
Then $\{v_i\}_{i=1}^{2m}$ and $\{v_{j}^{\ast}\}_{j=1}^{2m}$ are dual
bases for $V$ in the sense that $\langle v_i,
v_j^{\ast}\rangle=\delta_{i,j}$ for any $i,j$.\smallskip

There is a right action of $\mbb_n$ on $V^{\otimes n}$ which is
defined on generators by $$
\begin{aligned} (v_{i_1}\otimes\cdots\otimes v_{i_n})s_j
&:=-(v_{i_1}\otimes\cdots\otimes v_{i_{j-1}}\otimes
v_{i_{j+1}}\otimes v_{i_{j}}\otimes v_{i_{j+2}}
\otimes\cdots\\
&\qquad\qquad\otimes v_{i_n}),\\
(v_{i_1}\otimes\cdots\otimes v_{i_n})e_j
&:=-v_{i_1}\otimes\cdots\otimes v_{i_{j-1}}\otimes
\biggl(\sum_{k=1}^{2m}v_{k}\otimes v_k^{\ast}\biggr)\otimes
v_{i_{j+2}}\otimes\cdots\\ &\qquad\qquad \otimes
v_{i_n}.\end{aligned}
$$

Assume $n\geq 2$. For any pair of integers $1\leq s<t\leq n$ we
define the $(s,t)$-contraction operator $C_{s,t}: V^{\otimes
n}\rightarrow V^{\otimes n-2}$ by
$$\begin{aligned} (w_1\otimes\cdots\otimes
w_n)C_{s,t}&=\langle w_s,w_t\rangle w_1\otimes\cdots\otimes
w_{s-1}\otimes\widehat{w_s}\otimes w_{s+1}\otimes\cdots\otimes\\
&\qquad\qquad w_{t-1}\otimes\widehat{w_t}\otimes
w_{t+1}\otimes\cdots\otimes w_{n},\end{aligned}
$$
where $w_1,\cdots,w_n\in V$, $\widehat{w_s}, \widehat{w_t}$ mean
that we omit the tensor factors $w_s$ and $w_t$ in the tensor
product, and we define the $(s,t)$-expansion operator $D_{s,t}:
V^{\otimes n-2}\rightarrow V^{\otimes n}$ by
$$\begin{aligned} (w_1\otimes\cdots\otimes
w_{n-2})D_{s,t}&=\sum_{k=1}^{2m}w_1\otimes\cdots\otimes
w_{s-1}\otimes{v_k}\otimes w_{s}\otimes\cdots\otimes\\
&\qquad\qquad w_{t-1}\otimes{v_k^{\ast}}\otimes
w_{t}\otimes\cdots\otimes w_{n-2}.\end{aligned}
$$

\begin{lemma} {\rm (\cite[\S10.1.1]{GW})} For each pair of integers $1\leq s<t\leq n$, both $C_{s,t}$ and
$D_{s,t}$ are $Sp(V)$-module homomorphisms.
\end{lemma}

For each pair of integers $1\leq s<t\leq n$, we use $e_{s,t}$ to
denote the unique Brauer $n$-diagram which satisfies the following
two conditions:\begin{enumerate}
\item for any integer $a\in\{1,2,\cdots,n\}\setminus\{s,t\}$, the
vertex labeled by $a$ in the top row is connected with the vertex
labeled by $a$ in the bottom row;
\item the vertex labeled by $s$ in either the top row or the bottom row is connected with the vertex
labeled by $t$ in the same row.
\end{enumerate}
In particular, $e_s=e_{s,s+1}$ for any integer $1\leq s<n$.

\begin{lemma} \label{contra} For any pair of integers $1\leq s<t\leq n$, we have $\Ker e_{s,t}=\Ker
C_{s,t}$. In particular, $$ \bigcap_{1\leq s<t\leq n}\Ker
e_{s,t}=\bigcap_{1\leq s<t\leq n}\Ker C_{s,t}.
$$
\end{lemma}
\begin{proof} By definition, for any integers $1\leq i_1,\cdots,i_n\leq 2m$,
$$
(v_{i_1}\otimes\cdots\otimes
v_{i_n})e_{s,t}=-(v_{i_1}\otimes\cdots\otimes
v_{i_n})C_{s,t}D_{s,t}.
$$
It follows that for any $v\in V^{\otimes n}$, $v
e_{s,t}=-vC_{s,t}D_{s,t}$. On the other hand, it is clear that
$D_{s,t}$ maps the basis of simple $(n-2)$-tensors to a set of
$K$-linearly independent elements in $V^{\otimes n}$. Hence it is an
injective map from $V^{\otimes n-2}$ to $V^{\otimes n}$. It follows
that $\Ker e_{s,t}=\Ker C_{s,t}$. In particular,
$$ \bigcap_{1\leq s<t\leq n}\Ker e_{s,t}=\bigcap_{1\leq s<t\leq
n}\Ker C_{s,t}. $$
\end{proof}

\begin{lemma} \label{br} Let $x\in\BBd_n$. The following three statements are
equivalent:\begin{enumerate}
\item $x$ contains exactly two horizontal edges (one
edge in each of the top and the bottom rows in the diagram);
\item $x=ye_{s,t}$, for some $y\in\bBS_n$ and two integers $1\leq
s<t\leq n$;
\item $x=e_{s',t'}z$ for some $z\in\bBS_n$ and two integers $1\leq
s'<t'\leq n$.
\end{enumerate}
\end{lemma}
\begin{proof} By \cite{E}, $x$ contains exactly two horizontal edges
if and only if $$x=d_1^{-1}e_1\sigma d_2$$ for some $d_1,
d_2\in\mathfrak{D}_{(2,n-2)}$, $\sigma\in\bBS_{\{3,4,,\cdots,n\}}$,
where $\mathfrak{D}_{(2,n-2)}$ is the set of distinguished right
coset representatives of $\bBS_{\{1,2\}}\times\bBS_{\{3,4,\cdots,n\}}$ in
$\bBS_n$. It follows that
$$ x=\bigl(d_1^{-1}\sigma d_2\bigr)\bigl(d_2^{-1}e_1d_2\bigr)=
\bigl(d_1^{-1}e_1d_1\bigr)\bigl(d_1^{-1}\sigma d_2\bigr).
$$
Since both $d_2^{-1}e_1d_2$ and $d_1^{-1}e_1d_1$ are of the form
$e_{s,t}$ for two distinct integers $1\leq s<t\leq n$, the lemma
follows at once.
\end{proof}

\begin{definition} For each integer $1\leq f\leq
[n/2]$, we denote by $\mbb_n^{(f)}$ the two-sided ideal of $\mbb_n$
generated by $e_1e_3\cdots e_{2f-1}$.
\end{definition}
Note that $\mbb_n^{(f)}$ is spanned by all the Brauer $n$-diagrams
which contain at least $2f$ horizontal edges. Recall the definition of $W_{1,n}$ given in Definition \ref{w1n}.

\begin{corollary} \label{2dfn} With the notations as above, we have that $$
W_{1,n}=\bigcap_{1\leq s<t\leq n}\Ker C_{s,t}.
$$
\end{corollary}
\begin{proof} This follows directly from Lemma \ref{contra} and
\ref{br}.
\end{proof}

Note that the above corollary ensures that our definition \ref{w1n}
of harmonic tensors coincides with that given in \cite{DS} and
\cite[\S 10.2.1]{GW}.\smallskip

Let ``$\ast$" be the $K$-algebra anti-automorphism (of order two) of $\mbb_n$ which is
defined on generators by $$ s_i^{\ast}=s_i,\quad
e_i^{\ast}=e_i,\quad \forall\, 1\leq i\leq n-1.
$$
For any right $\mbb_n$-module $M$, the dual space
$M^{\ast}:=\Hom_{K}(M,K)$ is naturally endowed with a right
$\mbb_n$-module structure via the anti-involution ``$\ast$". That is,
$(fx)(v):=f(vx^{\ast}), \forall\,f\in M^{\ast}, x\in\mbb_n, v\in M$.
For each integer $1\leq k\leq 2m$, we set
$$I(k,n):=\bigl\{(i_1,\cdots,i_n)\bigm|i_j\in\{1,2,\cdots,k\},\,\,\forall\,
j\bigr\}.$$ For any $\mui=(i_1,\cdots,i_n)\in I(k,n)$, we write
$v_{\mui}=v_{i_1}\otimes\cdots\otimes v_{i_n}$. The bilinear form
$\langle\,\, ,\,\,\rangle$ on $V$ naturally induces a non-degenerate
bilinear form on $V^{\otimes n}$ such that $$ \langle v_{\mui},
v_{\muj}\rangle:=\prod_{s=1}^{n}\langle v_{i_s},
v_{j_s}\rangle,\quad\,\,\forall\, \mui,\muj\in I(2m,n).
$$
By definition and an easy check, we see that for any $v,w\in V^{\otimes n}, g\in
Sp(V), x\in\mbb_n$,
\begin{equation}\label{asso} \langle
gv,w\rangle=\langle v, g^{-1}w\rangle,\qquad \langle
vx,w\rangle=\langle v, w x^{\ast}\rangle .
\end{equation}
In other words, the bilinear form $\langle,\,\rangle$ induces a
$Sp(V)$-$\mbb_n$-bimodule isomorphism $$\begin{aligned}
\Theta:\,\,\,&
V^{\otimes n}\cong \Bigl(V^{\otimes n}\Bigr)^{\ast}\\
& v_{\mui}\mapsto \Theta(v_{\mui}):\,\,v_{\muj}\mapsto\langle
v_{\mui}, v_{\muj}\rangle,\,\,\forall\,\mui,\muj\in I(2m,n).
\end{aligned}
$$
By definition and (\ref{asso}), we deduce that $\Theta$ restricts to
an inclusion
\begin{equation}\label{embed}\Theta_1:\,\,
W_{1,n}\hookrightarrow \Bigl(V^{\otimes n}/V^{\otimes
n}\mbb_n^{(1)}\Bigr)^{\ast}.
\end{equation}
In the next section, we shall show that $\Theta_1$ is actually an isomorphism, see Corollary \ref{comdim}.

\section{Proof of Theorem \ref{mainthm0}, \ref{mainthm1} and Corollary \ref{maincor1}}\smallskip

The purpose of this section is to give a proof of Theorem
\ref{mainthm0}, \ref{mainthm1} and Corollary \ref{maincor1}. One of
the key steps is the use of Lemma \ref{Hu2lem}, whose proof will be
given in the next section.
\medskip

Let $G:=Sp_{2m}({K})$. We identify $V$ with $K^{2m}$ and $Sp(V)$
with $G$ by using the ordered basis $\bigl(v_1,v_2,\cdots,v_{2m}\bigr)$. Then $$
G=\Bigl\{A\in GL_{2m}(K)\Bigm|A^{t}JA=J\Bigr\},
$$
where $A^t$ denotes the matrix transpose of $A$, $$
J=\sum_{i=1}^{m}E_{i,2m+1-i}-\sum_{i=m+1}^{2m}E_{i,2m+1-i},
$$
and for each $1\leq i,j\leq 2m$, $E_{i,j}$ denotes the matrix unit which is $1$ on the $(i,j)$th position and $0$
elsewhere.
Recall that $G$ is a connected semisimple linear algebraic group over ${K}$. Let $K^{\times}:=K\setminus\{0\}$.
Let $T$ be the subgroup consisting of all diagonal matrices in $G$. That is, $$
T=\Bigl\{\sum_{i=1}^{m}\bigl(t_iE_{i,i}+t_{i}^{-1}E_{2m+1-i,2m+1-i}\bigr)\Bigm|t_1,\cdots,t_m\in K^{\times}\Bigr\}.
$$
Then $T$ is a maximal torus of $G$. Let $X(T)=\Hom_{G}(T,K^{\times})$ be the weight lattice of $G$. An element in $X(T)$
will also be called a $T$-weight. Let $X^{+}=X(T)^{+}$ be the set of dominant $T$-weights of $G$. Let $T_{2m}$ be the maximal torus consisting of all the diagonal matrices in $GL_{2m}(K)$. For each $1\leq i\leq 2m$,  let $\varepsilon_i$ be the function which sends a matrix in $T_{2m}$ to its $i$th element in the diagonal. Then $T$ is a subtorus of $T_{2m}$. Let $S:=\{\eps_i-\eps_{i+1}, 2\eps_m|1\leq i\leq m-1\}$, which is the set of simple roots of $Sp_{2m}(K)$. Let $B$ be the corresponding positive Borel subgroup. We refer the reader to \cite[15.2, Page 144--145]{DM} for the explicit description of $B$. For each $\lambda\in X^+$, we denote by $L(\lambda), \nabla(\lambda)$ and $\Delta(\lambda)$ the simple
module, co-Weyl module and Weyl module for $G$ associated to $\lam$ respectively. We identify the weight $\lam=\lam_1\eps_1+\cdots+\lam_m\eps_m\in X(T)$ with $(\lam_1,\cdots,\lam_m)\in\mmZ^m$. Thus a weight $\sum_{i=1}^{m}\lam_i\eps_i$ is dominant if and only if $(\lam_1,\lam_2,\cdots,\lam_m)$ is a partition.
For any $\lam,\mu\in X^{+}$, we
define
$$ \lam\geq\mu\Longleftrightarrow\lam-\mu\in\sum_{\alpha\in
S}\mathbb{Z}^{\geq 0}\alpha.
$$

For each simple $n$-tensor $v_{\mui}=v_{i_1}\otimes \cdots\otimes
v_{i_n}\in V^{\otimes n}$ with $\mui=(i_1,\cdots,i_n)\in I(2m,n)$,
the $T$-weight of $v_{\mui}$ is $$ \sum_{j=1}^n\sum_{1\leq i_j\leq
m}\eps_{i_j}-\sum_{j=1}^n\sum_{m+1\leq i_j\leq 2m}\eps_{2m+1-i_j}.
$$

Let $f$ be an integer with $0\leq f\leq [n/2]$. Let $\pi_{f}$ be the set of
dominant $T$-weights of $G$ appearing in $V^{\otimes n-2f}$. By \cite{W}, \begin{equation}\label{pif}
\pi_f=\biggl\{\lam_1\varepsilon_1+\cdots+\lam_m\varepsilon_m\biggm|\begin{matrix}
\text{$\lam=(\lam_1,\cdots,\lam_m)\vdash n-2f-2t$, for}\\
\text{some integer $0\leq t\leq [n/2-f]$.}
\end{matrix}\biggr\}.
\end{equation}
It is well-known that $\pi_f$ is saturated in the sense of \cite[Part II,
A.2]{Ja}. Let $\mathcal{C}(\pi_f)$ denote the category of
finite dimensional $G$-modules $M$ such that all composition factors
of $M$ have the form $L(\mu)$ with $\mu\in\pi_f$. Following
\cite[Part II, Chapter A]{Ja}, we define a functor
$\mathcal{O}_{\pi_f}$ from the category of finite dimensional
$G$-modules to the category $\mathcal{C}(\pi_f)$ as follows: $$
\mathcal{O}_{\pi_f}(M)=\sum_{M'\subseteq M, M'\in
\mathcal{C}(\pi_f)}M'.
$$

Let $M$ be a finite dimensional $G$-module. Recall that (\cite[Part
II, 4.16]{Ja}) an ascending chain
$$
0=M_0\subset M_1\subset M_2\subset\cdots\subset M_{s-1}\subset M_s=M
$$
of $G$-submodules is called a good filtration of $M$ if each
$M_i/M_{i-1}$ is isomorphic to some $\nabla(\lam_i)$ with $\lam_i\in
X^{+}$. We use $\bigl(M:\nabla(\lam_i)\bigr)$ to denote the number
of factors in the filtration isomorphic to $\nabla(\lam_i)$.
Similarly, an ascending chain
$$
0=N_0\subset N_1\subset N_2\subset\cdots\subset N_{t-1}\subset N_t=M
$$
of $G$-submodules is called a Weyl filtration of $M$ if each
$N_i/N_{i-1}$ is isomorphic to some $\Delta(\mu_i)$ with $\mu_i\in
X^{+}$. We use $\bigl(M:\Delta(\mu_i)\bigr)$ to denote the number of
factors in the filtration isomorphic to $\Delta(\mu_i)$.

\begin{lemma} \label{Jakey} {\rm (\cite[Part II, Chapter A]{Ja})} 1) A finite dimensional $G$-module $M$ belongs to
$\mathcal{C}(\pi_f)$ if and only if each dominant weight of $M$
belongs to $\pi_f$;

2) If $M$ is a finite dimensional $G$-module with a good filtration,
then $\mathcal{O}_{\pi_f}(M)$ has a good filtration with $$
\bigl(\mathcal{O}_{\pi_f}(M):\nabla(\lam)\bigr)=\begin{cases}
\bigl(M:\nabla(\lam)\bigr), &\text{if $\lam\in\pi_{f}$;}\\
0, &\text{otherwise;}
\end{cases}
$$

3) $\mathcal{O}_{\pi_f}$ is a left exact functor.
\end{lemma}

\begin{definition} \label{alpha} We set $\alpha:=\sum_{k=1}^{2m}v_k\otimes v_{k}^{\ast}\in V^{\otimes
2}$.
\end{definition}

Recall that $G$ acts diagonally on $V^{\otimes 2}$.

\begin{lemma} \label{Jakey2} 1) For any $g\in G$, we have that $g\alpha=\alpha$. That is, $K\alpha$ is a trivial $G$-module;

2) $\mathcal{O}_{\pi_f}\bigl(V^{\otimes
n}\mbb_n^{(f)}\bigr)=V^{\otimes n}\mbb_n^{(f)}$.
\end{lemma}

\begin{proof} As a $G$-module, $V^{\otimes 2}\cong V\otimes V^{\ast}\cong\End_K(V)$, where $\alpha$ was mapped to the identity map on $V$. It follows
that $G$ acts trivially on $K\alpha$. This proves 1).

Note that the left action of $G$ on $V^{\otimes n}$ commutes with
the right action of $\mbb_n$. By definition,  $V^{\otimes
n}\mbb_n^{(f)}$ is a sum of some submodules of the form $$
\bigl((K\alpha)^{\otimes f}\otimes V^{\otimes n-2f}\bigr)\sigma,
$$
where $\sigma\in\bBS_{n}$. Applying the result 1) that we have
proved, we get that (as a $G$-module)
$$ \bigl((K\alpha)^{\otimes f}\otimes V^{\otimes
n-2f}\bigr)\sigma\cong (K\alpha)^{\otimes f}\otimes V^{\otimes
n-2f}\cong V^{\otimes n-2f}.
$$
Therefore, it follows from definition that
$\mathcal{O}_{\pi_f}\bigl(V^{\otimes n}\mbb_n^{(f)}\bigr)=V^{\otimes
n}\mbb_n^{(f)}$.
\end{proof}

We want to show that $\mathcal{O}_{\pi_f}\bigl(V^{\otimes
n}\bigr)=V^{\otimes n}\mbb_n^{(f)}$. Before proving this equality, we
need some preparation. For simplicity, for any two finite
dimensional $G$-modules $M, N$, we use $\sum_{\phi: M\rightarrow
N}\Image\phi$ to denote the sum of all the image subspaces
$\Image\phi$, where $\phi$ runs over all the $G$-module
homomorphisms from $M$ to $N$. It is clearly a $G$-submodule of $N$.
We use $\tau_f$ to denote the following embedding:
\begin{equation}\label{tauf}
\begin{aligned}
\tau_f:\,\,\,V^{\otimes n-2f}&\hookrightarrow V^{\otimes n},\\
v_{i_1}\otimes\cdots\otimes v_{i_{n-2f}}&\mapsto
\alpha^{\otimes f}\otimes v_{i_1}\otimes\cdots\otimes v_{i_{n-2f}}.\end{aligned}\end{equation}

\begin{lemma} \label{keylem1} The $K$-linear space
$\Hom_{KG}\bigl(V^{\otimes n-2f}, V^{\otimes n}\mbb_n^{(f)}\bigr)$ is
spanned by all $\sigma{\tau_f}$, where $\sigma\in\bBS_n$. Furthermore,
the natural embedding $$ \iota_1: \Hom_{KG}\bigl(V^{\otimes n-2f},
V^{\otimes n}\mbb_n^{(f)}\bigr)\rightarrow\Hom_{KG}\bigl(V^{\otimes
n-2f}, V^{\otimes n}\bigr)
$$
is actually an isomorphism, and $$ \sum_{\phi: V^{\otimes
n-2f}\rightarrow V^{\otimes n}}\Image\phi=V^{\otimes
n}\mbb_n^{(f)}=\sum_{\phi: V^{\otimes n-2f}\rightarrow V^{\otimes
n}\mbb_n^{(f)}}\Image\phi.$$
\end{lemma}

\begin{proof} By Lemma \ref{Jakey2}, it is clear that $$
\sigma{\tau_f}\in\Hom_{KG}\bigl(V^{\otimes n-2f}, V^{\otimes
n}\mbb_n^{(f)}\bigr)
$$
for all $\sigma\in\bBS_n$.

By the proof of Lemma \ref{Jakey2}, we know that $V^{\otimes
n}\mbb_n^{(f)}$ is a sum of some submodules, each of which is
isomorphic to $V^{\otimes n-2f}$. It follows that $$ V^{\otimes
n}\mbb_n^{(f)}=\sum_{\phi: V^{\otimes n-2f}\rightarrow V^{\otimes
n}\mbb_n^{(f)}}\Image\phi\subseteq \sum_{\phi: V^{\otimes
n-2f}\rightarrow V^{\otimes n}}\Image\phi.
$$

As a $G$-module, $V\cong V^{\ast}$. It follows that
$$\begin{aligned} \Hom_{KG}\bigl(V^{\otimes n-2f}, V^{\otimes
n}\bigr) &\cong \Bigl(\bigl(V^{\otimes n-2f}\bigr)^{\ast}\otimes
V^{\otimes
n}\Bigr)^{G}\\
&\cong\bigl(V^{\otimes n-2f}\otimes V^{\otimes
n}\bigr)^{G}=\bigl(V^{\otimes
2n-2f}\bigr)^{G}\cong\End_{KG}\bigl(V^{\otimes
n-f}\bigr).\end{aligned}
$$
Here the isomorphisms in the above equalities are induced from the
following natural isomorphisms: $$\begin{aligned}
&\rho:\,\,V^{\otimes
2n-2f}\overset{\sim}{\rightarrow}\End_{K}\bigl(V^{\otimes
n-f}\bigr)\\
&\quad v_{\mui}\mapsto\rho(v_{\mui}):\,\,v_{j_1}\otimes\cdots\otimes
v_{j_{n-f}}\mapsto\Bigl(\prod_{s=1}^{n-f}\langle
v_{i_s},v_{j_s}\rangle\Bigr)v_{i_{n-f+1}}\otimes\cdots\otimes
v_{i_{2n-2f}};\\
&\rho_f:\,\,V^{\otimes
2n-2f}\overset{\sim}{\rightarrow}\Hom_{K}\bigl(V^{\otimes
n-2f},V^{\otimes n}\bigr)\\
&\quad v_{\mui}\mapsto\rho(v_{\mui}):\,\,v_{j_1}\otimes\cdots\otimes
v_{j_{n-2f}}\mapsto\Bigl(\prod_{s=1}^{n-2f}\langle
v_{i_s},v_{j_s}\rangle\Bigr)v_{i_{n-2f+1}}\otimes\cdots\otimes
v_{i_{2n-2f}},
\end{aligned}
$$
where $v_{\mui}:=v_{i_1}\otimes\cdots\otimes
v_{i_{2n-2f}}$.\smallskip

Recall that there is a natural right action of the Brauer algebra
$\mbb_{n-f}$ on the $(n-f)$-tensor space $V^{\otimes n-f}$. We use
$\varphi'_K$ to denote the natural $K$-algebra homomorphism from
$(\mbb_{n-f})^{\rm op}$ to $\End_{KG}\bigl(V^{\otimes n-f}\bigr)$. By
\cite[Theorem 1.4]{DDH}, we know that $\varphi'_K$ is surjective.
Therefore, $\End_{KG}\bigl(V^{\otimes n-f}\bigr)$ is spanned by the
elements $\varphi'_K(D)$, where $D$ runs over the Brauer
$(n-f)$-diagrams in $\BBd_{n-f}$.

Let $D\in\BBd_{n-f}$. We regard $\varphi'_K(D)$ as an element in
$V^{\otimes 2n-2f}$ using the isomorphism $\rho$. By the
definition of $\varphi'_K$, it is easy to see that
$$ \varphi'_K(D)=\pm \bigl(\alpha^{\otimes n-f}\bigr)\sigma,
$$
for some $\sigma\in\bBS_{2n-2f}$. Now using the isomorphism $\rho_f$ one can check the image of $\varphi'_K(D)$.
It follows that $\Hom_{KG}\bigl(V^{\otimes n-2f}, V^{\otimes n}\bigr)$ is
spanned by all $\sigma{\tau_f}$, where $\sigma\in\bBS_n$.

It remains to prove that $\Hom_{KG}\bigl(V^{\otimes n-2f},
V^{\otimes n}\mbb_n^{(f)}\bigr)=\Hom_{KG}\bigl(V^{\otimes
n-2f}, V^{\otimes n}\bigr)$. It suffices to show that the image of $\rho_f\bigl(\varphi'_K(D)\bigr)$ is
contained in $V^{\otimes n}\mbb_n^{(f)}$ for each $D\in\BBd_{n-f}$. Note that $\sum_{k=1}^{2m}v_{k}\otimes
v_{k}^{\ast}=-\sum_{k=1}^{2m}v_{k}^{\ast}\otimes v_{k}$. In
particular, we can write
$$ \varphi'_K(D)=\pm\sum_{1\leq i_1,\cdots,i_{n-f}\leq
2m}\cdots\otimes v_{i_j}\otimes\cdots\otimes
v_{i_j}^{\ast}\otimes\cdots\in V^{\otimes 2n-2f},
$$
where $i_1,\cdots,i_{n-f}$ are $n-f$ independent summation indices.
The positions of each pair of $(v_{i_j}, v_{i_j}^{\ast})$ in the
above sum are uniquely determined by $\sigma$ and hence by $D$. For
example, if \medskip
\begin{center}
\scalebox{0.25}[0.25]{\includegraphics{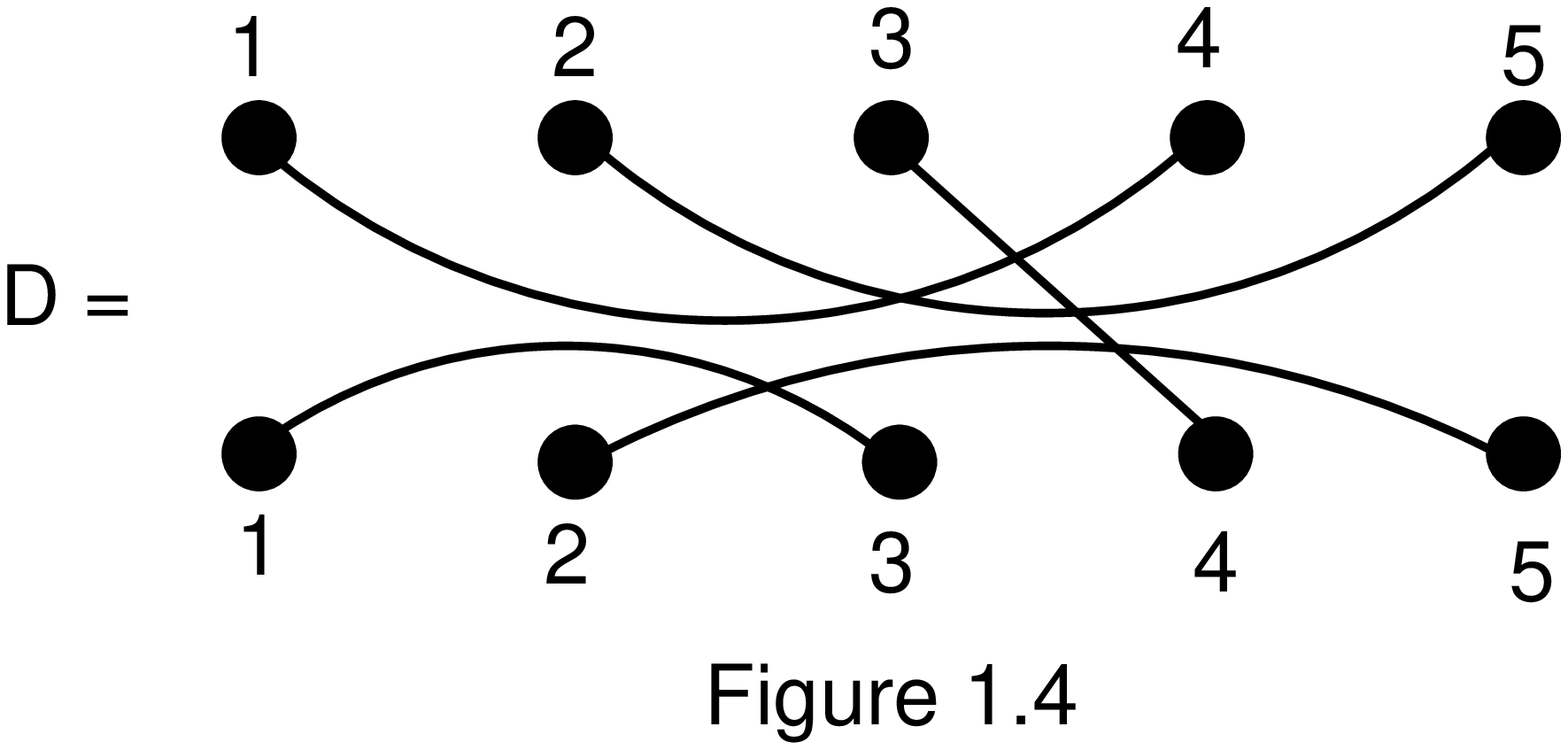}}
\end{center}
\medskip
then $$ \varphi'_K(D)=-\sum_{1\leq i_1,i_2,i_3,i_4,i_{5}\leq
2m}v_{i_1}\otimes v_{i_2}\otimes v_{i_3}\otimes
v_{i_1}^{\ast}\otimes v_{i_2}^{\ast}\otimes v_{i_4}\otimes
v_{i_5}\otimes v_{i_4}^{\ast}\otimes v_{i_3}^{\ast}\otimes
v_{i_5}^{\ast}.
$$
Each independent summation index $i_j$ was attached with a pair of
integers $1\leq a_j<b_j\leq 2n-2f$ such that $v_{i_j},
v_{i_j}^{\ast}$ appear in the position $a_j, b_j$ respectively. Using \cite[Lemma 4.2]{Hu1}
and the definition of $\rho_f$, it suffices to show that
$$\#\Bigl\{(a_j,b_j)\Bigm|a_{j}\geq n-2f+1\Bigr\}\geq f. $$

We fix a bijection between the horizontal edges in the top row of
$D$ and the horizontal edges in the bottom row of $D$. We consider a
vertex $A$ in the top row of $D$ which is labeled by an integer $x$
with $n-2f+1\leq x\leq n-f$. There are only three possibilities:
\begin{enumerate}
\item $A$ is connected with a vertex $B$ in the top row of $D$ which is labeled by an integer $y$ with $n-2f+1\leq
y\leq n-f$, then the horizontal edge $(A,B)$ must correspond to a horizontal edge $(A',B')$ in
the bottom row of $D$. Applying \cite[Lemma 4.2]{Hu1} and the
definition of $\rho$, we see that the edges $(A,B), (A',B')$ determine two integers $j, k$, such that $n-2f+1\leq a_j<b_j\leq n-f,\,\,n-f+1\leq a_k<b_k\leq 2n-2f$;
\item $A$ is connected with a vertex $B$ in the top row of $D$ which is labeled by an integer $y$ with $1\leq
y\leq n-2f$, then the horizontal edge $(A,B)$ must correspond to a horizontal edge $(A',B')$ in
the bottom row of $D$. Applying \cite[Lemma 4.2]{Hu1} and the
definition of $\rho$, we see that the edges $(A',B')$ determines an integer $j$,
such that $n-f+1\leq a_j<b_j\leq 2n-2f$;
\item $A$ is connected with a vertex $B$ in the bottom row of $D$, then applying \cite[Lemma 4.2]{Hu1} and the
definition of $\rho$, we see that the vertical edge $(A,B)$ determines an integer $j$,
such that $n-f+1\leq a_j=x\leq n-f<b_j\leq 2n-2f$.
\end{enumerate}
By a simple counting, we deduce that $\#\Bigl\{(a_j,b_j)\Bigm|a_{j}\geq
n-2f+1\Bigr\}\geq f$. This completes the proof of the lemma.
\end{proof}

Let $U$ be the unipotent radical of $B$. For any $G$-module $M$, a vector $v\in M$
is said to be a maximal vector if $xv=v$ for any $x\in U$.

\begin{lemma} \label{BP} {\rm (\cite[4.13 (3), B.3 (5), B.4 Remark]{Ja})} Let $\lam,\mu\in X^{+}$, $i\in\mathbb{Z}^{\geq 0}$.
\begin{enumerate}
\item $\Ext_{G}^{i}\bigl(\Delta(\lam),\nabla(\mu)\bigr)=\delta_{i,0}\delta_{\lam,\mu}K$;
\item $\Ext_{G}^{1}\bigl(\nabla(\lam),
\nabla(\mu)\bigr)\neq 0$ implies that $\lam>\mu$. In particular,
$$\Ext_{G}^{1}\bigl(\nabla(\lam), \nabla(\lam)\bigr)=0;$$
\item $\Ext_{G}^{1}\bigl(\Delta(\lam),
\Delta(\mu)\bigr)\neq 0$ implies that $\lam<\mu$. In particular,
$$\Ext_{G}^{1}\bigl(\Delta(\lam), \Delta(\lam)\bigr)=0;$$
\item $\Hom_{KG}\bigl(\nabla(\lam),\nabla(\mu)\bigr)\neq 0$ unless $\mu\leq\lam$, and
$\Hom_{KG}\bigl(\nabla(\lam),\nabla(\lam)\bigr)=K$.
\end{enumerate}
\end{lemma}

\begin{lemma} \label{BP2} Let $M$ be a $G$-module and $\lam\in X^{+}$. Suppose $\rho: \nabla(\lam)\rightarrow
M$ is a surjective $G$-module homomorphism such that
$\rho(u_{\lam})\neq 0$ for some highest weight vector
$u_{\lam}\in\nabla(\lam)$ of weight $\lam$, then $\rho$ is an isomorphism. In
particular, $M\cong\nabla(\lam)$.
\end{lemma}

\begin{proof} This follows from the fact that
$\soc_{G}\bigl(\nabla(\lam)\bigr)=L(\lam)$.
\end{proof}

Recall that for each element $\lambda=(\lam_1,\cdots,\lam_m)\in\mmZ^m$, $\lam$ is identified with the weight
$\lam_1\varepsilon_1+\cdots+\lam_m\varepsilon_m\in X(T)$.

\begin{lemma} \label{compare} Let $a, b$ be two integers such that $0\leq a<b\leq [n/2]$. Let $\lam$ be a partition of
$n-2a$ and $\mu$ be a partition of $n-2b$. Suppose that $\ell(\lam),
\ell(\mu)\leq m$. Then $\lam\nleqslant\mu$.
\end{lemma}

\begin{proof} Suppose that $\lam\leq\mu$. Then there exist some integers $a_1,\cdots,a_m\geq 0$ such that $$
\sum_{i=1}^m \mu_i\varepsilon_i-\sum_{i=1}^m \lam_i\varepsilon_i=\sum_{i=1}^{m-1}a_i\varepsilon_i+2a_m\varepsilon_m.
$$
Comparing the coefficients of $\varepsilon_i$ for each $1\leq i\leq m$ and adding all of them together, we get that
$$
(n-2b)-(n-2a)=\sum_{i=1}^m\mu_i-\sum_{i=1}^m\lam_i=\sum_{i=1}^{m-1}a_i+2a_m\geq 0,
$$
which is a contradiction to the assumption that $a<b$.
\end{proof}

As a $G$-module, $V=\Delta(\eps_1)=\nabla(\eps_1)=L(\eps_1)$ is a tilting module. By \cite[Proposition E.7]{Ja}, we know that
$V^{\otimes n}$ is a tilting module too. In particular,
$V^{\otimes n}$ has both a Weyl filtration and a good filtration. We
fix an integer $f$ with $0\leq f\leq [n/2]$. In view of Lemma
\ref{BP} and Lemma \ref{compare}, we can find a good filtration of
$V^{\otimes n}$:
\begin{equation} \label{filtra1} 0=M_0\subset M_1\subset M_2\subset\cdots\subset
M_p=V^{\otimes n},
\end{equation}
such that \begin{enumerate}
\item $M_i/M_{i-1}\cong \underbrace{\nabla(\lam^{(i)})\oplus\cdots\oplus\nabla(\lam^{(i)})}_{\text{$n_i$ copies}}$
for $i=1,2,\cdots,p$, where $n_i\in\mathbb{N}$, $\lam^{(i)}\vdash
n-2g_i$ for some integers $0\leq g_i\leq [n/2]$ and
$\ell(\lam^{(i)})\leq m$;
\item $\lam^{(i)}\neq\lam^{(j)}$ for any $i\neq j$, and $\lam^{(i)}<\lam^{(j)}$ only if
$i<j$;
\item there exists an integer $1\leq t\leq p$,
such that $\lam^{(i)}\in \pi_f$ if and only if $1\leq i\leq t$.
\end{enumerate}
We refer the reader to \cite[Part II, 4.16, Remark 4]{Ja}) for the construction of the filtration.

\begin{lemma} \label{keylem2} We keep the notations as above. Let $i$ be an integer with
$1\leq i\leq p$.

1) If $v$ is a (non-zero) maximal vector of weight $\lam^{(i)}$ in
$V^{\otimes n}$, then $v+M_{i-1}$ is also a (non-zero) maximal vector
of weight $\lam^{(i)}$ in $V^{\otimes n}/M_{i-1}$. Moreover, every
maximal vector of weight $\lam^{(i)}$ in $V^{\otimes n}/M_{i-1}$
arises in this way;

2) if $N$ is a $G$-submodule of $V^{\otimes n}/M_{i-1}$ which is
isomorphic to $\nabla(\lam^{(i)})$, then $N\subseteq M_{i}/M_{i-1}$;

3) $M_ib\subseteq M_i$ for any $b\in\mbb_n$. In other words, the
filtration (\ref{filtra1}) is actually a $KG$-$\mbb_n$-bimodules
filtration.
\end{lemma}

\begin{proof} 1) This is a combination of \cite[Part II, Proposition 4.16 (a), Lemma 2.13 (a)]{Ja}.

2) Since $V^{\otimes n}/M_{i}$ has a good filtration with each
section being isomorphic to some $\nabla(\mu)$ with
$\mu\nleqq\lam^{(i)}$, it follows from Lemma \ref{BP} that $$
\Hom_{KG}\bigl(N, V^{\otimes n}/M_{i}\bigr)=0.
$$
Therefore, the natural embedding $$\Hom_{KG}\bigl(N,
M_{i}/M_{i-1}\bigr)\hookrightarrow\Hom_{KG}\bigl(N, V^{\otimes
n}/M_{i-1}\bigr)
$$
becomes an equality, from which we deduce that $N\subseteq
M_{i}/M_{i-1}$.

3) This follows from the following simple fact: If $M$ is a finite dimensional $G$-module with a good filtration and
$\lam\in X^{+}$ is minimal such that
$\nabla(\lam)$ appears in a good filtration of $M$, then $$
M':=\sum_{\phi\in\Hom_{KG}(\nabla(\lam),M)}\phi(\nabla(\lam))
$$
is a direct sum of copies of $\nabla(\lam)$ and $M/M'$ has a good filtration without factor $\nabla(\lam)$. In particular, $M'$ is a right module for $\End_{KG}(M)$ and $M/M'$ has an $\End_{KG}(M)$-module structure.
\end{proof}

Let $k$ be an integer and $\lam$ a composition of $k$.
The conjugate of $\lam$ is the partition $\lam'=(\lam'_1,\lam'_2,\cdots)$ where
$\lam'_i=\#\{j\geq 1|\lam_j\geq i\}$. The Young diagram of $\lam$ is the set $$
[\lam]:=\bigl\{(a,b)\bigm|1\leq b\leq\lam_a\bigr\}.
$$
Assume that $\lam$ is a partition of $n$. A $\lam$-tableau is a bijective map $\mathfrak{t}: [\lam]\rightarrow\{1,2,\dots,n\}$. A standard $\lam$-tableau
is a $\lam$-tableau in which the entries increase along each row and down each column.
Let $\mathfrak{t}^{\lam}$ (resp., $\mathfrak{t}_{\lam}$) be the
standard $\lam$-tableau such that the numbers $1,2,\cdots,n$ appear
in order along the rows (resp., along the columns). Let $\bBS_{\lam}$ be the Young subgroup of $\bBS_n$ corresponding to
$\lam$, which
is the subgroup fixing the sets $\{1,2,\cdots,\lam_1\},\{\lam_1+1,\lam_1+2,\cdots,\lam_1+\lam_2\},\cdots$.
Let $\mathcal{D}_{\lam}$ be the set of distinguished right coset representatives of $\bBS_{\lam}$ in $\bBS_n$.
Let
$w_{\lam}\in\bBS_n$ such that
$\mathfrak{t}^{\lam}w_{\lam}=\mathfrak{t}_{\lam}$. Then $w_{\lam}\in\mathcal{D}_{\lam}$. Let
$x_{\lam}:=\sum_{w\in\bBS_{\lam}}w$. Recall the definition of $\alpha$ in Definition \ref{alpha}.

\begin{definition} \label{Zdfn} Let $g$ be an integer with $0\leq g\leq [n/2]$ and $\lam$ be a
partition of $n-2g$ with $\ell(\lam)\leq m$. We set
$$\begin{aligned} v_{\lam}:&=\underbrace{v_1\otimes\cdots\otimes
v_1}_{\text{$\lam_1$ copies}}\otimes\cdots\otimes
\underbrace{v_m\otimes\cdots\otimes
v_m}_{\text{$\lam_m$ copies}};\\
z_{g,\lam}:&=\alpha^{\otimes g}\otimes v_{\lam}w_{\lam}x_{\lam'}.
\end{aligned}
$$
\end{definition}

Recall that $U$ is the unipotent radical of the positive Borel subgroup $B$. For each $G$-module $M$, we use $M^{U}_{\lam}$ to denote the
subspace of maximal vectors in $M$ of weight $\lam$. The following
lemma plays a key role in this section. The proof will be given in
the next section.

\begin{lemma} \label{Hu2lem} Let $g$ be an integer with $0\leq g\leq [n/2]$ and $\lam$ be a
partition of $n-2g$ with $\ell(\lam)\leq m$. Then $\bigl(V^{\otimes
n}\bigr)^{U}_{\lam}=z_{g,\lam}\mbb_n$. In particular, $z_{g,\lam}$
is a non-zero maximal vector of weight $\lam$ in $V^{\otimes n}$.
Moreover, the dimension of $\bigl(V^{\otimes n}\bigr)^{U}_{\lam}$
(and hence of $z_{g,\lam}\mbb_n$) is independent of $K$;
\end{lemma}

The following proposition proves the first part of Theorem
\ref{mainthm3}.

\begin{proposition} \label{mainprop2} Let $i$ be an integer with $1\leq i\leq p$. Then $z_{g_i,\lam^{(i)}}\in M_i$
and there is a $KG$-$\mbb_n$-bimodule isomorphism: $$
M_i/M_{i-1}\cong \nabla(\lam^{(i)})\otimes z_{g_i,\lam^{(i)}}\mbb_n.
$$
\end{proposition}

\begin{proof} By construction, $M_i/M_{i-1}\cong\bigl(\nabla(\lam^{(i)}))^{\oplus n_i}$. We fix a decomposition: $$
M_i/M_{i-1}=N_{i,1}\oplus\cdots\oplus N_{i,n_i},$$ where $N_{i,s}\cong\nabla(\lam^{(i)})$ for each $1\leq s\leq n_i$.
Let $u_0$ be a fixed non-zero maximal vector of weight $\lam^{(i)}$ in $\nabla(\lam^{(i)})$. For each $1\leq s\leq n_i$, let $u_s$ be a fixed non-zero maximal vector of weight $\lam^{(i)}$ in $N_{i,s}$ and let $\eta_s: \nabla(\lam^{(i)})\rightarrow N_{i,s}$ be the unique isomorphism which sends $u_0$ to $u_s$.
Then it is easy to see (by definition and comparing dimensions) that the following map
$$ \begin{aligned}\widetilde{\psi}_i:\,\,\nabla(\lam^{(i)})\otimes\bigl(V^{\otimes n}\bigr)_{\lam^{(i)}}^{U} &\rightarrow
M_i/M_{i-1}\\
v\otimes\sum_{s=1}^{n_i}a_su_s&\mapsto
\sum_{s=1}^{n_i}a_s\eta_s(v),\quad\,\forall\,v\in\nabla(\lam^{(i)}),\,a_s\in K,
\end{aligned}
$$
is a left $KG$-module isomorphism. We claim that it is also a right $\mbb_n$-module homomorphism.

It suffices to show that $\widetilde{\psi}_i(v\otimes u_sb)=\eta_s(v)b$ for each $1\leq s\leq n_i$, $b\in\mbb_n$ and $v\in\nabla(\lam^{(i)})$.
We write $u_sb=\sum_{t=1}^{n_i}a_tu_t$, where $a_t\in K$ for each $t$. Then by definition, $\widetilde{\psi}_i(v\otimes u_sb)=\sum_{t=1}^{n_i}a_t\eta_t(v)$.
By the commuting action between $KG$ and $\mbb_n$, it is easy to check that $\sum_{t=1}^{n_i}a_t\eta_t-\eta_sb$ is a left $KG$-homomorphism from $\nabla(\lam^{(i)})$ to $M_i/M_{i-1}$. By direct verification, $\Bigl(\sum_{t=1}^{n_i}a_t\eta_t-\eta_sb\Bigr)(u_0)=0$, hence
$\Bigl(\sum_{t=1}^{n_i}a_t\eta_t-\eta_sb\Bigr)(v)=0$ for any $v\in KGu_0=L(\lam^{(i)})=\soc_G\bigl(\nabla(\lam^{(i)})\bigr)$.
Since $M_i/M_{i-1}\cong\bigl(\nabla(\lam^{(i)})\bigr)^{\oplus n_i}$ and $\End_{KG}\bigl(\nabla(\lam^{(i)})\bigr)=K$, it follows that any non-zero map in
$\Hom_{KG}\bigl(\nabla(\lam^{(i)}),M_i/M_{i-1}\bigr)$ must be injective. We deduce that
$\Bigl(\sum_{t=1}^{n_i}a_t\eta_t-\eta_sb\Bigr)(v)=0$ for any $v\in\nabla(\lam^{(i)})$, as required. This proves that $\widetilde{\psi}_i$
is a right $\mbb_n$-module homomorphism as well.

Finally, by Lemma \ref{Hu2lem}, $\bigl(V^{\otimes n}\bigr)^{U}_{\lam^{(i)}}=z_{g_i,\lam^{(i)}}\mbb_n$. Hence the proposition follows.
\end{proof}

\begin{lemma} Let $t$ be the integer which was introduced in the third property of the filtration (\ref{filtra1}). We have that $$
V^{\otimes n}\mbb_n^{(f)}\subseteq M_t.
$$
\end{lemma}
\begin{proof} By Lemma \ref{keylem1}, it is enough to show that $$
\sum_{\phi: V^{\otimes n-2f}\rightarrow V^{\otimes
n}}\Image\phi\subseteq M_t.
$$
Since $V^{\otimes n-2f}$ has a Weyl filtration such that each
section is isomorphic to some $\Delta(\mu)$ with $\mu\in\pi_f$ and
$V^{\otimes n}/M_t$ has a good filtration such that each section is
isomorphic to some $\nabla(\nu)$ with $\nu\notin\pi_f$, it follows
that $$ \Hom_{KG}\bigl(V^{\otimes n-2f}, V^{\otimes
n}/M_t\bigr)=0.
$$
Hence $$ \Hom_{KG}\bigl(\sum_{\phi: V^{\otimes
n-2f}\rightarrow V^{\otimes n}}\Image\phi, V^{\otimes
n}/M_t\bigr)=0,
$$
which implies that $$ \sum_{\phi: V^{\otimes n-2f}\rightarrow
V^{\otimes n}}\Image\phi\subseteq M_t,
$$
as required. \end{proof}

Our purpose is to show that $V^{\otimes n}\mbb_n^{(f)}=M_t$.

\begin{theorem} \label{keystep3} $M_t=\sum_{\phi: V^{\otimes n-2f}\rightarrow V^{\otimes
n}}\Image\phi=V^{\otimes
n}\mbb_n^{(f)}=\mathcal{O}_{\pi_f}\bigl(V^{\otimes n}\bigr)$. In
particular, both $\dim V^{\otimes n}\mbb_n^{(f)}$ and $\dim
V^{\otimes n}/V^{\otimes n}\mbb_n^{(f)}$ are independent of $K$.
\end{theorem}

\begin{proof} First, by Lemma \ref{Jakey} 2), we have that $$
\dim \mathcal{O}_{\pi_f}\bigl(V^{\otimes n}\bigr)=\sum_{\lam\in\pi_f}(V^{\otimes
n}:\nabla(\lam))\dim\nabla(\lam).
$$
Note that the character $\chf(\nabla(\lam))$ of $\nabla(\lam)$ is given by the Weyl character formula. In particular, $\dim\nabla(\lam)$ is independent of $K$. Note also that the character $\chf(V^{\otimes n})$ is equal to $(\chf(V))^n$ which is also independent of $K$. So the filtration multiplicities
$(V^{\otimes n}:\nabla(\lam))$ is nothing but the coefficient of the Weyl character $\chi_{\lam}$ in the expansion of $(\chf(V))^n$ into a linear combination of Weyl characters. In particular, $$
(V^{\otimes n}:\nabla(\lam))=[V_{\mmC}^{\otimes n}:\nabla(\lam)_{\mmC}],
$$
which is independent of $K$. Since the set $\pi_f$ is also independent of $K$, we conclude that $\dim \mathcal{O}_{\pi_f}\bigl(V^{\otimes n}\bigr)$ is independent of $K$. Moreover, by the construction of $M_t$, it is clear that $$
\dim M_t=\sum_{\lam\in\pi_f}(V^{\otimes n}:\nabla(\lam))\dim\nabla(\lam)=\dim \mathcal{O}_{\pi_f}\bigl(V^{\otimes n}\bigr).
$$

Now, by Lemma \ref{keylem1}, $\sum_{\phi: V^{\otimes n-2f}\rightarrow V^{\otimes n}}\Image\phi=V^{\otimes n}\mbb_n^{(f)}$. By Lemma \ref{Jakey2} 2) and
Lemma \ref{Jakey} 3), $V^{\otimes n}\mbb_n^{(f)}=\mathcal{O}_{\pi_f}\bigl(V^{\otimes n}\mbb_n^{(f)}\bigr)\subseteq\mathcal{O}_{\pi_f}\bigl(V^{\otimes n}\bigr)$. Therefore, to complete the proof of the theorem, it suffices to prove that
$$M_t\subseteq\sum_{\phi: V^{\otimes n-2f}\rightarrow V^{\otimes
n}}\Image\phi. $$

We use induction on $i$ to show that $M_i\subseteq \sum_{\phi:
V^{\otimes n-2f}\rightarrow V^{\otimes n}}\Image\phi$ for each
integer $1\leq i\leq t$. Suppose that $M_{i-1}\subseteq \sum_{\phi:
V^{\otimes n-2f}\rightarrow V^{\otimes n}}\Image\phi$.

Recall that $t$ is equal to the number of dominant weights in $\pi_f$. Using a similar construction of the filtration (\ref{filtra1}), we
can get a good filtration of $V^{\otimes n-2f}$:
\begin{equation}\label{filtra2}
0=M'_0\subseteq M'_1\subseteq\cdots\subseteq
M'_{t}=V^{\otimes n-2f}, \end{equation} such that
\begin{enumerate}
\item $M'_j/M'_{j-1}\cong \underbrace{\nabla(\lam^{(j)})\oplus\cdots\oplus\nabla(\lam^{(j)})}_{\text{$n'_j$ copies}}$
for $j=1,2,\cdots,t$, where $n'_j\in\mathbb{Z}^{\geq 0}$,
$\lam^{(j)}\vdash n-2f-2h_j$ for some integers $0\leq h_j\leq
[n/2-f]$ and $\ell(\lam^{(j)})\leq m$;
\item $\lam^{(i)}\neq\lam^{(j)}$ for any $i\neq j$, and $\lam^{(i)}<\lam^{(j)}$ only if
$i<j$.
\end{enumerate}
Let $j$ be an integer with $1\leq j\leq t$. By the properties of the
filtration (\ref{filtra1}) and the filtration (\ref{filtra2}) we
have constructed and Lemma \ref{BP} 4), it is clear that
$\Hom_{KG}\bigl(M'_j, V^{\otimes n}/M_j\bigr)=0$. It follows that the
natural embedding $$ \Hom_{KG}\bigl(M'_j,
M_j\bigr)\hookrightarrow\Hom_{KG}\bigl(M'_j, V^{\otimes n}\bigr)
$$
becomes an equality. By Lemma \ref{Jakey2}, we see that the following map  $$ \iota_f: V^{\otimes
n-2f}\rightarrow V^{\otimes n}, \quad v\mapsto\alpha^{\otimes
f}\otimes v,\,\,\forall\,v\in V^{\otimes n-2f},
$$ is an embedding. It follows that $\iota_f(M'_j)\subseteq M_j$
for each $1\leq j\leq t$. In particular, $\iota_f$ induces a natural
map from $M'_i/M'_{i-1}$ to $M_i/M_{i-1}$.

Recall that for each $1\leq i\leq p$, $\lam^{(i)}\vdash n-2g_i$. In particular,
$g_i=f+h_i$ for each $1\leq i\leq t$. For each $1\leq i\leq t$, we define $z'_{h_i,\lam^{(i)}}:=
\alpha^{\otimes h_i}\otimes v_{\lam}w_{\lam}x_{\lam'}$. Then $z'_{h_i,\lam^{(i)}}$ is a
maximal vector of weight $\lam^{(i)}$ in $V^{\otimes n-2f}$ and
$z'_{h_i,\lam^{(i)}}\in M'_i$ (cf. Lemma \ref{Hu2lem} and
Proposition \ref{mainprop2}). It is clear that
$\iota_f(z'_{h_i,\lam^{(i)}})=z_{g_i,\lam^{(i)}}$. Let $N'$ be a
$G$-submodule of $M'_i/M'_{i-1}$ such that $z'_{h_i,\lam^{(i)}}\in
N'\cong\nabla(\lam^{(i)})$. Applying Lemma \ref{BP2}, we deduce that
$\iota_f(N')\cong\nabla(\lam^{(i)})$. Therefore, combining with our
induction hypothesis, we can find a $G$-submodule $\widetilde{N}:=M_{i-1}+\iota_f(M'_i)$ of
$V^{\otimes n}$ such that $$ M_{i-1}\subseteq
\widetilde{N}\subseteq\sum_{\phi: V^{\otimes n-2f}\rightarrow
V^{\otimes n}}\Image\phi,\,\,z_{g_i,\lam^{(i)}}\in\widetilde{N},\,\,
N:=\widetilde{N}/M_{i-1}\cong\nabla(\lam^{(i)}).
$$
By the proof of Proposition \ref{mainprop2}, we can deduce that $M_i/M_{i-1}\subseteq N\mbb_n$.
As a consequence, we get that
$$ M_i\subseteq\sum_{\phi: V^{\otimes n}\rightarrow V^{\otimes
n}}\phi(\widetilde{N}).
$$
Since $\widetilde{N}\subseteq\sum_{\phi: V^{\otimes n-2f}\rightarrow
V^{\otimes n}}\Image\phi$, it follows that $M_i\subseteq\sum_{\phi:
V^{\otimes n-2f}\rightarrow V^{\otimes n}}\Image\phi$, as required.
\end{proof}

\begin{corollary}\label{maincor25} With the notations as before, we have
that $$ V^{\otimes n}=V^{\otimes n}\mbb_n^{(f)}+\sum_{\substack{0\leq
g<f\\ \mu\vdash n-2g}}(KG)z_{g,\mu}\mbb_n.
$$
\end{corollary}
\begin{proof} This follows directly from the proof of Theorem
\ref{keystep3}.
\end{proof}

\bigskip

With Theorem \ref{keystep3} in hand, Theorem \ref{mainthm0} and \ref{mainthm1} follow almost trivially from some standard arguments
in commutative algebras. For the reader's convenience, we include the details below.
\medskip\medskip

\noindent {\bf Proof of Theorem \ref{mainthm0}, \ref{mainthm1} and
Corollary \ref{maincor1}, \ref{maincor15}:}\smallskip

By Theorem \ref{keystep3} and Lemma \ref{Jakey}, we know that
$V^{\otimes n}\mbb_n^{(f)}$ has a good filtration. Note that
$V^{\otimes n}$ has a good filtration too. Applying \cite[Part II,
4.17]{Ja}, we deduce that both $V^{\otimes n}/V^{\otimes
n}\mbb_n^{(f)}$ and $V^{\otimes n}\mbb_n^{(f-1)}/V^{\otimes
n}\mbb_n^{(f)}$ have a good filtration. This proves Theorem
\ref{mainthm1} and the first statement of Corollary \ref{maincor15}.\medskip

It is clear that $V_{\mmZ}^{\otimes n}\mbb_n^{(f)}$ is a free
$\mmZ$-module of finite rank $r_f$. Since $V_{\mmZ}^{\otimes
n}\mbb_n^{(f)}\otimes_{\mmZ}\mmC\cong V_{\mmC}^{\otimes n}\mbb_n^{(f)}$, it
follows that $r_f=\dim V_{\mmC}^{\otimes n}\mbb_n^{(f)}$. By Theorem
\ref{keystep3}, we know that $\dim V^{\otimes n}\mbb_n^{(f)}$ is
independent of $K$ and hence $\dim V^{\otimes
n}\mbb_n^{(f)}=r_f$. Thus the canonical surjection $V_{\mmZ}^{\otimes
n}\mbb_n^{(f)}\otimes_{\mmZ}K\rightarrow V^{\otimes n}\mbb_n^{(f)},
v\otimes a\mapsto va$ must be an isomorphism.

Now let $F$ be an arbitrary field and $\overline{F}$ be the algebraic closure of $F$. Recall that $V_{F}:=V_{\mmZ}\otimes_{\mmZ}F$, $V_{\overline{F}}:=V_{\mmZ}\otimes_{\mmZ}\overline{F}$. It is clear that the canonical map
$V_{F}^{\otimes n}\mbb_n^{(f)}\otimes_{F}\overline{F}\rightarrow
V_{\overline{F}}^{\otimes n}\mbb_n^{(f)}$ is an isomorphism. It
follows that $\dim V_{F}^{\otimes n}\mbb_n^{(f)}=\dim
V_{\overline{F}}^{\otimes n}\mbb_n^{(f)}=r_f$. In particular, the
first statement of Corollary \ref{maincor1} also follows. For any
commutative $\mmZ$-algebra $F$ which is a field, the canonical
surjection $V_{\mmZ}^{\otimes n}\mbb_n^{(f)}\otimes_{\mmZ}F\rightarrow
V_{F}^{\otimes n}\mbb_n^{(f)}$ must be an isomorphism. For any
commutative $\mmZ$-algebra $R$, let $F$ be a field such that $F\cong
R/m$ for some maximal ideal $m$ of $R$. We have the following
commutative diagram of maps:
\begin{equation}\label{cd1}\begin{CD}
V_{\mmZ}^{\otimes n}\mbb_n^{(f)}\otimes_{\mmZ}R\otimes_{R}F
@>{\pi_1\otimes\newid}>>
V_{R}^{\otimes n}\mbb_n^{(f)}\otimes_{R} F\\
@V{\newid\otimes\pi} V{\wr} V @V{\pi_2}VV\\
V_{\mmZ}^{\otimes n}\mbb_n^{(f)}\otimes_{\mmZ}F @>{\sim}>{\pi_3}>
V_{F}^{\otimes n}\mbb_n^{(f)}
\end{CD},
\end{equation}
 where $\pi, \pi_1, \pi_2, \pi_3$ are all canonical maps. Since
$\pi_3$ is an isomorphism, it follows that the canonical surjection
$\pi_1: V_{\mmZ}^{\otimes n}\mbb_n^{(f)}\otimes_{\mmZ}R\rightarrow
V_{R}^{\otimes n}\mbb_n^{(f)}$ must be injective and hence an isomorphism. This proves the first isomorphism in
part 2) of Theorem \ref{mainthm0}.\medskip

Using the same argument as before, we know that $$\begin{aligned}
\dim V_{F}^{\otimes n}/V_{F}^{\otimes n}\mbb_n^{(f)}&=\dim
V_{\overline{F}}^{\otimes n}/V_{\overline{F}}^{\otimes
n}\mbb_n^{(f)}\\
&=\dim V_{\mmC}^{\otimes n}/V_{\mmC}^{\otimes n}\mbb_n^{(f)}=(2m)^n-r_f.
\end{aligned}$$ Note that
$\rank V_{\mmZ}^{\otimes n}/V_{\mmZ}^{\otimes n}\mbb_n^{(f)}=\dim
V_{\mmQ}^{\otimes n}/V_{\mmQ}^{\otimes n}\mbb_n^{(f)}=(2m)^n-r_{f}$.
In order to show that $V_{\mmZ}^{\otimes n}/V_{\mmZ}^{\otimes
n}\mbb_n^{(f)}$ is a free $\mmZ$-module of rank $(2m)^n-r_f$, we
consider the following commutative diagram of maps:
\begin{equation}\label{cd2}
\xymatrix{V_{\mmZ}^{\otimes n}\mbb_n^{(f)}\otimes_{\mmZ}F
\ar[r]^{\quad\iota\otimes\newid} \ar[d]^{\wr} & V_{\mmZ}^{\otimes
n}\otimes_{\mmZ}F
\ar@{>>}[r]^{\pi\otimes\newid\quad\quad\,}\ar[d]^{\wr} &
V_{\mmZ}^{\otimes n}/V_{\mmZ}^{\otimes
n}\mbb_n^{(f)}\otimes_{\mmZ}F \ar@{>>}[d]^{\theta}\\
V_{F}^{\otimes n}\mbb_n^{(f)} \ar@{^{(}->}[r] & V_{F}^{\otimes n}
\ar@{>>}[r] & V_{F}^{\otimes n}/V_{F}^{\otimes n}\mbb_n^{(f)},}
\end{equation}
where $\iota$ denotes the natural injection $V_{\mmZ}^{\otimes
n}\mbb_n^{(f)}\hookrightarrow V_{\mmZ}^{\otimes n}$, $\pi$ denotes the
natural projection $V_{\mmZ}^{\otimes n}\twoheadrightarrow
V_{\mmZ}^{\otimes n}/V_{\mmZ}^{\otimes n}\mbb_n^{(f)}$. By diagram
chasing, it is easy to see that the natural surjection $\theta:
V_{\mmZ}^{\otimes n}/V_{\mmZ}^{\otimes
n}\mbb_n^{(f)}\otimes_{\mmZ}F\twoheadrightarrow V_{F}^{\otimes
n}/V_{F}^{\otimes n}\mbb_n^{(f)}$ is an injection and hence an
isomorphism. In particular, $$ \dim V_{\mmZ}^{\otimes
n}/V_{\mmZ}^{\otimes n}\mbb_n^{(f)}\otimes_{\mmZ}F=\dim V_{F}^{\otimes
n}/V_{F}^{\otimes n}\mbb_n^{(f)}=(2m)^n-r_f.
$$
We claim that $V_{\mmZ}^{\otimes n}/V_{\mmZ}^{\otimes n}\mbb_n^{(f)}$
must be a free $\mmZ$-module. In fact, suppose this is not the case,
then $V_{\mmZ}^{\otimes n}/V_{\mmZ}^{\otimes n}\mbb_n^{(f)}$ must contain
a non-zero $p$-torsion element for some prime number $p$. It follows
that $$ \dim_{\mathbb{F}_p} \Bigl(V_{\mmZ}^{\otimes n}/V_{\mmZ}^{\otimes
n}\mbb_n^{(f)}\otimes_{\mmZ}\mathbb{F}_p\Bigr)>(2m)^n-r_f,
$$
which is a contradiction. This proves our claim. As a consequence,
we conclude that $V_{\mmZ}^{\otimes n}\mbb_n^{(f)}$ is a pure
$\mmZ$-submodule of $V_{\mmZ}^{\otimes n}$.

Now using a commutative diagram similar to (\ref{cd1})
(replacing $V_{\mmZ}^{\otimes n}\mbb_n^{(f)}$ with $V_{\mmZ}^{\otimes
n}/V_{\mmZ}^{\otimes n}\mbb_n^{(f)}$), we can argue as before that for
any commutative $\mmZ$-algebra $R$, the canonical map $V_{\mmZ}^{\otimes
n}/V_{\mmZ}^{\otimes n}\mbb_n^{(f)}\otimes_{\mmZ}R\rightarrow
V_{R}^{\otimes n}/V_{R}^{\otimes n}\mbb_n^{(f)}$ is always an
isomorphism. This completes the proof of Theorem \ref{mainthm0}.\medskip

It remains to prove the second statement of Corollary
\ref{maincor1}. Let $K$ be an algebraically closed field. Recall the bimodule isomorphism $\Theta$ introduced in
the paragraph below (\ref{asso}). $\Theta$ induces a map $\Theta_f$
from $\mathcal{HT}_f^{\otimes n}$ to $\bigl(V^{\otimes
n}\mbb_n^{(f)}\bigr)^{\ast}$. The second equality in (\ref{asso}) and
the definition of $\mathcal{HT}_f^{\otimes n}$ imply that the image
of $\Theta_f$ is contained in $\bigl(V^{\otimes
n}\mbb_n^{(f)}/V^{\otimes n}\mbb_n^{(f+1)}\bigr)^{\ast}$. Therefore,
$\Theta_f$ is a bimodule homomorphism from $\mathcal{HT}_f^{\otimes
n}$ to $\bigl(V^{\otimes n}\mbb_n^{(f)}/V^{\otimes
n}\mbb_n^{(f+1)}\bigr)^{\ast}$. We claim that $\Theta_f$ is
injective.

Let $u\in \mathcal{HT}_f^{\otimes
n}$ such that $\Theta_f(u)=0$. Now for any $v\in
V^{\otimes n}$, we can write (by applying Corollary \ref{maincor25})
$$ v=v_0+\sum_{\substack{h\in G,\,\, c\in\mbb_n\\ 0\leq g<f,\,\, \mu\vdash
n-2g}}hz_{g,\mu}c,
$$
where $v_0\in V^{\otimes n}\mbb_n^{(f)}$. Using the assumption that $\Theta_f(u)=0$ and applying the two equalities in
(\ref{asso}) we get that $$
\langle u, v\rangle=\langle u, v_0\rangle+\sum_{\substack{h\in G,\,\, c\in\mbb_n\\ 0\leq g<f,\,\, \mu\vdash
n-2g}}\langle h^{-1}uc^{\ast},z_{g,\mu}\rangle=0+0=0.
$$
Since the bilinear form $\langle\, ,\rangle$ on $V^{\otimes n}$ is non-degenerate, it follows that
$u=0$. Hence $\Theta_f$ is injective. This proves our claim.

To emphasize the base field, we use the notations
$\mathcal{HT}_{f,K}^{\otimes n}$, $\mathcal{HT}_{f,\mmC}^{\otimes n}$.
Note that we have already known that the dimension of $V^{\otimes
n}\mbb_n^{(f)}$ is independent of the field $K$. The space
$\mathcal{HT}_f^{\otimes n}$ can be identified as the solution space
of a homogeneous system of linear equations whose coefficients are
all defined over $\mmZ$. In particular, $$
\dim\mathcal{HT}_{f,\mmC}^{\otimes n}\leq\dim\mathcal{HT}_{f,K}^{\otimes n}.
$$
By \cite[(10.3.14)]{GW}, we know that $$
\dim\mathcal{HT}_{f,\mmC}^{\otimes n}=\dim V_{\mmC}^{\otimes
n}\mbb_n^{(f)}/V_{\mmC}^{\otimes n}\mbb_n^{(f+1)}=\dim V^{\otimes
n}\mbb_n^{(f)}/V^{\otimes n}\mbb_n^{(f+1)}.
$$
By the injectivity of $\Theta_f$, we know that $$
\dim\mathcal{HT}_{f,K}^{\otimes n}\leq\dim V^{\otimes
n}\mbb_n^{(f)}/V^{\otimes n}\mbb_n^{(f+1)}.
$$
It follows that $\dim\mathcal{HT}_{f,K}^{\otimes n}=\dim V^{\otimes
n}\mbb_n^{(f)}/V^{\otimes n}\mbb_n^{(f+1)}$. Therefore, $\Theta_f$
must be an isomorphism. This completes the proof of Corollary
\ref{maincor1}.
\medskip

Recall the definition of $\Theta_1$ in (\ref{embed}). The following result is a consequence of Corollary \ref{maincor1}.

\begin{corollary} \label{comdim} With the notations as above, we have that $$
\dim W_{1,n}^{\mmC}=\dim W_{1,n}=\dim V^{\otimes n}/V^{\otimes
n}\mbb_n^{(1)}.
$$
In particular,  $\Theta_1$ is always an isomorphism.
\end{corollary}
\begin{proof} Note that the space $W_{1,n}$ can be identified as the solution
space of a homogeneous system of linear equations whose coefficients
are all defined over $\mmZ$. This implies that $\dim W_{1,n}^{\mmC}\leq \dim W_{1,n}$.
Now applying Corollary \ref{maincor1} and (\ref{embed}), we prove the corollary.
\end{proof}

\section{Proof of Lemma \ref{Hu2lem} and Theorem \ref{mainthm3} 2)}

The purpose of this section is to give a proof of Lemma \ref{Hu2lem}
as well as Theorem \ref{mainthm3} 2). The key step in our proof is
to show that (for any $\lam\vdash n$)
$(KG)z_{0,\lam}\cong\Delta(\lam)$ and $V^{\otimes n}/(KG)z_{0,\lam}$
has a Weyl filtration as a $KG$-module. The proof makes use of
Lusztig's theory of canonical basis and based modules. To this end,
we have to work in a quantized setting.
\smallskip

Let $q$ be an indeterminate over $\mmZ$. Let $\mmA:=\mmZ[q,q^{-1}]$ be the
ring of Laurent polynomials in $q$. Let $\mbb_n(-q^{2m+1},q)_{\mmA}$ be
the specialized Birman-Murakami-Wenzl algebra (\cite{BW},
\cite{M}). By definition, it has generators $T_1,\cdots,T_{n-1},
E_1,\cdots,E_{n-1}$ which satisfy the following relations:
\begin{enumerate}
\item $T_i-T_i^{-1}=(q-q^{-1})(1-E_i)$, for $1\le i\le n-1$,
\item $E_i^2=(1-\sum_{i=-m}^{m}q^{2i})E_i$, for $1\le i\le n-1$,
\item $T_{i}T_{i+1}T_i=T_{i+1}T_iT_{i+1}$,  for $1\le i\le n-2$,
\item $T_iT_j=T_jT_i$,  for $|i-j|>1$,
\item $E_iE_{i+1}E_i=E_i,\,\,E_{i+1}E_{i}E_{i+1}=E_{i+1}$,\,\, for $1\le i\le n-2$,
\item $T_iT_{i+1}E_{i}=E_{i+1}E_i,\,\,T_{i+1}T_{i}E_{i+1}=E_{i}E_{i+1}$,\,\,for $1\le i\le n-2$,
\item $ E_i T_i=T_iE_i=-q^{-2m-1}E_i$,  for $1\le i\le n-1$.
\item $E_{i}T_{i+1}E_{i}=-q^{2m+1}E_{i},\,\,E_{i+1}T_{i}E_{i+1}=-q^{2m+1}E_{i+1}$,\,\,for $1\le i\le n-2$.
\end{enumerate}
Let ${\mmU}_{\mmQ(q)}(\mathfrak{sp}_{2m})$ be the quantized
enveloping algebra of $\mathfrak{sp}_{2m}(\mmC)$ over $\mmQ(q)$. Let
$e_i,f_i,k_i,k_i^{-1}, 1\leq i\leq m$ be the Chevalley generators of
$\mmU_{\mmQ(q)}(\mathfrak{sp}_{2m})$. Then it is a Hopf algebra with
coproduct $\Delta$, counit $\varepsilon$ and antipode $S$ defined on
generators by
$$\begin{aligned} &\Delta(e_i)=e_i\otimes 1+\widetilde{k}_{i}\otimes e_i,\quad
\Delta(f_i)=1\otimes f_i+f_i\otimes
\widetilde{k}_{i}^{-1},\quad \Delta(k_i)=k_i\otimes k_i,\\
& \varepsilon(e_i)=\varepsilon(f_i)=0,\quad \varepsilon(k_i)=1,\\
& S(e_i)=-\widetilde{k}_{i}^{-1}e_i,\quad S(f_i)=-f_i\widetilde{k}_{i},\quad S(k_i)=k_i^{-1},
\end{aligned}$$
where $$
\widetilde{k}_i:=\begin{cases} k_i, & \text{if $i\neq m$;}\\
k_i^2, &\text{if $i=m$.}\end{cases}
$$
We regard $\mmZ$ as an $\mmA$-algebra by specializing $q$ to $1$. For
each $\mmA$-module $N$ and $x\in N$, we set $N_{\mmZ}:=N_{\mmA}\otimes_{\mmA}\mmZ$ and
$x\!\downarrow_{q=1}:=x\otimes_{\mmA}1_{\mmZ}\in N_{\mmZ}$. We call $x\!\downarrow_{q=1}$ the specialization of $x$ at $q=1$.
Let $\mathbb{U}_{\mmA}(\mathfrak{sp}_{2m})$ be Lusztig's $\mmA$-form in
$\mathbb{U}_{\mmQ(q)}(\mathfrak{sp}_{2m})$. Let $V_{\mmA}$ be the free
$\mmA$-module spanned by $\{v_i\}_{i=1}^{2m}$. It is well-known that
there are natural commuting actions (cf. \cite[Section 2,3]{Hu2}) of
$\mathbb{U}_{\mmA}(\mathfrak{sp}_{2m})$ and $\mbb_n(-q^{2m+1},q)_{\mmA}$
on $V_{\mmA}^{\otimes n}$. If we specialize $q$ to $1$, then each
$T_i$ (resp., $E_i$) specializes to $-s_i$ (resp., to $e_i$), and
$\mbb_n(-q^{2m+1},q)$ will become the Brauer algebra $\mbb_n(-2m)$ in
this paper. Moreover, the action of $\mbb_n(-q^{2m+1},q)$ on
$V_{\mmA}^{\otimes n}$ becomes the action of $\mbb_n(-2m)$ on
$V_{\mmZ}^{\otimes n}$ in this paper. We recall that the representation
$\varphi_{C}$ of $\mbb_n(-q^{2m+1},q)_{\mmA}$ on $V_{\mmA}^{\otimes n}$ is
defined on generators as follows:
$$ \varphi_C(T_j):=\Bid_{V^{\otimes
j-1}}\otimes\beta'\otimes\Bid_{V^{\otimes n-j-1}},\quad
\varphi_C(E_j):=\Bid_{V^{\otimes
j-1}}\otimes\gamma'\otimes\Bid_{V^{\otimes n-j-1}},
$$
for all $1\leq j\leq n-1$, where
$$\begin{aligned} &\beta':=\sum_{1\leq i\leq
2m}\Bigl(qE_{i,i}\otimes E_{i,i}+q^{-1}E_{i,i'}\otimes
E_{i',i}\Bigr)+\sum_{\substack{1\leq
i,j\leq 2m\\ i\neq j,j'}} E_{i,j}\otimes E_{j,i}+\\
&\qquad\qquad (q-q^{-1})\sum_{1\leq i<j\leq 2m}\Bigl(E_{i,i}\otimes
E_{j,j}-q^{\rho_j-\rho_i}\epsilon_i\epsilon_j E_{i,j'}\otimes
E_{i',j}\Bigr),\\
&\gamma':=\sum_{1\leq i,j\leq
2m}q^{\rho_j-\rho_i}\epsilon_i\epsilon_j E_{i,j'}\otimes E_{i',j},\\
&(\rho_1,\cdots,\rho_{2m}):=(m,m-1,\cdots,1,-1,\cdots,-m+1,-m),
\quad \epsilon_i:=\sign(\rho_i),
\end{aligned}
$$
and each $E_{i,j}$ is the matrix unit (i.e.,
$E_{i,j}v_{k}=\delta_{k,j}v_i$ for each $1\leq k\leq 2m$). Let $\mHH_{\mmA}(\bBS_n)$
be the Iwahori--Hecke algebra associated to the symmetric group
$\bBS_n$, defined over $\mmA$ and with parameter $q$. By definition,
$\mHH_{\mmA}(\bBS_n)$ has generators $\hat{T}_1,\cdots,\hat{T}_{n-1}$
which satisfy the following relations:
$$\begin{aligned} &(\hat{T}_i-q)(\hat{T}_i+q^{-1})=0,
\,\,\text{for $
i=1,2,\cdots,n-1$};\\
&
\hat{T}_i\hat{T}_{i+1}\hat{T}_i=\hat{T}_{i+1}\hat{T}_i\hat{T}_{i+1},\,\,
\text{for $1\leq i\leq n-2$;}\\
&\hat{T}_i\hat{T}_j=\hat{T}_j\hat{T}_i,\,\,\text{if $|i-j|>1$.}
\end{aligned}
$$
Let $\mHH_{q}(\bBS_n):=\mHH_{\mmA}(\bBS_n)\otimes_{\mmA}\mmQ(q)$.
For each composition $\lam$ of $n$, we use $\mHH_q(\bBS_{\lam})$ to
denote the Hecke algebra over $\mmQ(q)$ associated to the Young
subgroup $\bBS_{\lam}$. For each $w\in\bBS_n$, we use $\ell(w)$ to
denote the minimal integer $k$ such that $w=s_{j_1}s_{j_2}\cdots
s_{j_k}$; in that case, we call $s_{j_1}s_{j_2}\cdots s_{j_k}$ a
reduced expression of $w$. We define
$$
\hat{T}_w=\hat{T}_{j_1}\hat{T}_{j_2}\cdots \hat{T}_{j_k}\in
\mHH_{\mmA}(\bBS_n), \,\,T_w=T_{j_1}T_{j_2}\cdots T_{j_k}\in\mbb_n(-q^{2m+1},q)_{\mmA}, $$
if $w=s_{j_1}s_{j_2}\cdots s_{j_k}$ with $\ell(w)=k$. It is well-known that this is well defined, i.e.,
independent of the choice of the reduced expression. Let
$\hat{V}_{\mmA}$ be the free $\mmA$-submodule of $V_{\mmA}$ generated by
$\{v_i\}_{1\leq i\leq m}$. We recall that the representation
$\varphi_{A}$ of $\mHH_{\mmA}(\bBS_n)$ on $\hat{V}_{\mmA}^{\otimes n}$ is
defined on generators as follows: $$
\varphi_A(\hat{T}_j)=\Bid_{V^{\otimes
j-1}}\otimes\hat{\beta}\otimes\Bid_{V^{\otimes n-j-1}},\quad
\text{for $j=1,2,\cdots,n-1$,}
$$ where $$
\hat{\beta}:=\sum_{1\leq i\leq m}\Bigl(qE_{i,i}\otimes
E_{i,i}\Bigr)+\sum_{\substack{1\leq i,j\leq m\\ i\neq j}}
\Bigl(E_{i,j}\otimes E_{j,i}\Bigr)+(q-q^{-1})\sum_{1\leq i<j\leq
m}\Bigl(E_{i,i}\otimes E_{j,j}\Bigr).\\
$$
Let $\lam$ be a partition with $\ell(\lam)\leq m$. Recall the definition of $v_{\lam}$ in Definition \ref{Zdfn}.
We have that $v_{\lam}\hat{T}_{\sigma}=q^{\ell(\sigma)}v_{\lam}$ for
each $\sigma\in\bBS_{\lam}$.

\begin{definition} For each partition $\lam$ of $n$, we define
$$\begin{aligned} &X_{\lam}=\sum_{w\in\bBS_{\lam}}q^{\ell(w)}{T}_w,\quad
Y_{\lam'}=\sum_{w\in\bBS_{\lam'}}(-q)^{-\ell(w)}{T}_w.\\
&\hat{X}_{\lam}=\sum_{w\in\bBS_{\lam}}q^{\ell(w)}\hat{T}_w,\quad
\hat{Y}_{\lam'}=\sum_{w\in\bBS_{\lam'}}(-q)^{-\ell(w)}\hat{T}_w.
\end{aligned}$$
\end{definition}

\begin{definition} {\rm (\cite[Definition 3.1]{HX})} Set $$\alpha_q:=\sum_{1\leq k\leq 2m}q^{-\rho_k}\epsilon_k v_k\otimes
v_{k'}. $$ For each integer $0\leq g\leq [n/2]$ and each partition $\lam$ of
$n-2g$ with $\ell(\lam)\leq m$, we define $$
Z_{g,\lam}=\alpha_q^{\otimes g}\otimes v_{\lam}{T}_{w_{\lam}}Y_{\lam'}.
$$
\end{definition}

Recall that for each $w\in\bBS_n$, $T_w\!\downarrow_{q=1}=(-1)^{\ell(w)}w$ (instead of $w$). Thus
$z_{g,\lam}=Z_{g,\lam}\!\downarrow_{q=1}$. We set $V_{\mmQ(q)}:=V_{\mmA}\otimes_{\mmA}\mmQ(q)$
and $\hat{V}_{\mmQ(q)}:=\hat{V}_{\mmA}\otimes_{\mmA}\mmQ(q)$. Since $v_{\lam}\in\hat{V}_{\mmA}$, it is easy to check that $
v_{\lam}{T}_{w_{\lam}}Y_{\lam'}=v_{\lam}\hat{T}_{w_{\lam}}\hat{Y}_{\lam'}$.

\begin{lemma} {\rm (\cite[Lemmas 3.1, 3.2]{HX})} \label{lmhx} 1) $\mmQ(q)\alpha_q$ is the one dimensional trivial
$\mmU_{\mmQ(q)}(\mathfrak{sp}_{2m})$-submodule of
$V_{\mmQ(q)}^{\otimes 2}$, i.e., $x\alpha_q=\varepsilon(x)\alpha_q$
for any $x\in\mmU_{\mmQ(q)}(\mathfrak{sp}_{2m})$;

2) $Z_{g,\lam}$ is a non-zero maximal vector in $V_{\mmQ(q)}^{\otimes n}$ of
weight $\lam$ with respect to the action of ${\mmU}_{\mmQ(q)}(\mathfrak{sp}_{2m})$.
\end{lemma}
\begin{proof} This follows from direct verification and the fact
that the right action of $\mbb_n(-q^{2m+1},q)$ on $V_{\mmQ(q)}^{\otimes
n}$ commutes with the left action of
${\mmU}_{\mmQ(q)}(\mathfrak{sp}_{2m})$.
\end{proof}

 Let $\hat{M}:=(\hat{V}_{\mmQ(q)})^{\otimes {n}}$.
Let ${\mmU}_{\mmQ(q)}(\mathfrak{gl}_{m})$ be the quantized
enveloping algebra of the general linear Lie algebra
$\mathfrak{gl}_{m}(\mmC)$ over $\mmQ(q)$. Let
$\mathbb{U}_{\mmA}(\mathfrak{gl}_{m})$ be Lusztig's $\mmA$-form in
$\mathbb{U}_{\mmQ(q)}(\mathfrak{gl}_{m})$. There is a natural
representation of $\mathbb{U}_{\mmA}(\mathfrak{gl}_{m})$ on
$\hat{V}_{\mmA}$ and hence on $\hat{V}_{\mmA}^{\otimes n}$ (cf.
\cite{DPS}). By \cite[(27.3)]{Lu}, the
$\mmU_{\mmQ(q)}({\mathfrak{gl}_{m}})$-module $\hat{M}$ is a based
module. There is a bar involution $\psi_A$ which is defined on
$\hat{M}$. For any integers $1\leq i_1,\cdots,i_n\leq m$, there is a
unique element
${v}_{i_1}{\diamond}{v}_{i_2}{\diamond}\cdots{\diamond}
{v}_{i_{n}}\in \hat{V}_{\mmA}^{\otimes n}$, such that
\begin{enumerate}
\item $\psi_A\bigl({v}_{i_1}{\diamond}{v}_{i_2}{\diamond}\cdots{\diamond} {v}_{i_{n}}\bigr)={v}_{i_1}{\diamond}{v}_{i_2}{\diamond}\cdots{\diamond} {v}_{i_{n}}$, and
\item  ${v}_{i_1}{\diamond}\cdots{\diamond}{v}_{i_{n}}$ is equal to ${v}_{i_1}\otimes\cdots\otimes {v}_{i_{n}}$ plus a linear
combination of elements ${v}_{j_1}\otimes\cdots\otimes {v}_{j_{n}}$ with $({v}_{j_1},\cdots,{v}_{j_{n}})<^{A}({v}_{i_1},\cdots,{v}_{i_{n}})$
and with coefficients in $q^{-1}\mmZ[q^{-1}]$, where $``<^{A}"$ is a partial order\footnote{Here we use the subscript and superscript ``A'' to emphasize that both the bar involution and the partial order depend on the type $A_{m-1}$ root system.} defined in \cite[(27.3.1)]{Lu}.
\end{enumerate}

Let $\hat{B}:=\bigl\{{v}_{i_1}{\diamond}\cdots{\diamond}{v}_{i_{n}}\bigr\}_{1\leq i_1,\cdots,i_n\leq m}$. Then $\hat{B}$ is the set of Lusztig's canonical bases of $\hat{M}$. In particular, $\hat{B}$ is an $\mmA$-basis of $\hat{V}_{\mmA}^{\otimes n}$. Similarly,
${M}:=({V}_{\mmQ(q)})^{\otimes {n}}$ is a based module as a
${\mmU}_{\mmQ(q)}(\mathfrak{sp}_{2m})$-module. Let $B$ be the set of
canonical bases of $M$. Then $B$ is actually an $\mmA$-basis of
$M_{\mmA}:=V_{\mmA}^{\otimes n}$. Recall that $\pi_0$ is the set of
dominant $T$-weights appeared in $M$ (cf. (\ref{pif})). By \cite[(27.2.1)]{Lu}, there
is a partition $$ B=\bigsqcup_{\lam\in\pi_0}B[\lam].
$$
For each $\lam\in\pi_0$, we set $$
B[\geq\!\lam]:=\bigsqcup_{\lam\leq\mu\in\pi_0}B[\mu],\quad B[>\!\lam]:=\bigsqcup_{\lam<\mu\in\pi_0}B[\mu].
$$
Let $M[\geq\!\lam]_{\mmA}, M[>\!\lam]_{\mmA}$ be the $\mmA$-submodule of
$M_{\mmA}$ generated by the canonical basis elements in
$B[\geq\!\lam], B[>\!\lam]$ respectively. By \cite[27.1.8]{Lu}, both
$M[\geq\!\lam]_{\mmA}$ and $M[>\!\lam]_{\mmA}$ are
${\mmU}_{\mmA}(\mathfrak{sp}_{2m})$-stable and we have that
\begin{equation}\label{filtra3}
M[\geq\!\lam]_{\mmA}/M[>\!\lam]_{\mmA}\cong
\bigl(\Delta(\lam)_{\mmA}\bigr)^{\oplus n_{\lam}},
\end{equation}
for some $n_{\lam}>0$. Moreover, the canonical image of each $b\in
B[\lam]$ in $M[\geq\!\lam]_{\mmA}/M[>\!\lam]_{\mmA}$ is mapped to a
canonical basis element of some direct summand $\Delta(\lam)_{\mmA}$
in the right-hand side of (\ref{filtra3}).

Recall that a weight of ${\mmU}_{\mmQ(q)}(\mathfrak{gl}_{m})$ is
identified as an element $\mu=(\mu_1,\cdots,\mu_{m})\in\mmZ^{m}$ by
setting $\langle\mu,\eps_i\rangle= \mu_i$ for $1\leq i\leq m$; in this case, $\mu$ is said to be dominant if
$\mu_1\geq\mu_2\geq\cdots\geq\mu_m$. For any two weights $\lam, \mu$
of ${\mmU}_{\mmQ(q)}(\mathfrak{gl}_{m})$, we write $\lam>_{A}\mu$ if
$\sum_{i=1}^m(\lam_i-\mu_i)\eps_i=\sum_{i=1}^{m-1}a_i(\eps_i-\eps_{i+1})$
for some $a_1,\cdots,a_{m-1}\in\mmZ^{\geq 0}$. Let $\hat{\pi}_{0}$ be the set of dominant weights of
${\mmU}_{\mmQ(q)}(\mathfrak{gl}_{m})$ appeared in
$\hat{V}_{\mmQ(q)}^{\otimes n}$. Then $\hat{\pi}_0$ is nothing but the set of partitions of $n$ with no more than $m$ parts.
For each $\lam\in\hat{\pi}_{0}$, let
$\hat{B}[\lam]^{hi}$ be as defined in \cite[27.2.3]{Lu}. Using $\hat{\pi}_0$ and the partial order $>_{A}$, one can define the subsets
$\hat{B}[\geq\!\lam], \hat{B}[>\!\lam]$ and the $\mmA$-submodule $\hat{M}[\geq\!\lam]_{\mmA}, \hat{M}[>\!\lam]_{\mmA}$ in a similar way as before.
We set $$
\hat{B}^{hi}:=\bigsqcup_{\lam\in\hat{\pi}_{0}}\hat{B}[\lam]^{hi},\quad
{B}^{hi}:=\bigsqcup_{\lam\in {\pi}_{0}}{B}[\lam]^{hi}.
$$
Let $\mu=(\mu_1,\cdots,\mu_l)$ be a partition of $n$ satisfying $\mu_1\leq m$ and $\mu_l\neq 0$. We define $$
\widehat{b}_{\mu}:=\underbrace{v_{\mu_1}\diamond\cdots\diamond v_2\diamond v_1}_{\text{$\mu_1$ terms}}
\diamond\underbrace{v_{\mu_{2}}\diamond\cdots\diamond v_2\diamond v_1}_{\text{$\mu_{2}$ terms}}\diamond\cdots
\diamond\underbrace{v_{\mu_l}\diamond\cdots\diamond v_2\diamond v_1}_{\text{$\mu_l$ terms}}\in
\hat{V}_{\mmA}^{\otimes n},$$
which is (by definition) a canonical basis element in $\hat{B}$.

\begin{lemma} \label{lmkey41} With the notations as above, we have that $\widehat{b}_{\mu}\in\hat{B}^{hi}$.
\end{lemma}

\begin{proof} For each integer $1\leq i\leq m-1$, let $\widetilde{E}_i, \widetilde{F}_i$ denote the Kashiwara
operators as defined in
\cite{Lu}. Since we are following the notations and definitions of coproduct in Lusztig's book, the actions of
Kashiwara operators on the tensor product of crystal bases are given by the following formulae (cf.
\cite[17.2.4, 20.2.2, 27.3.3]{Lu})
$$\begin{aligned}
\widetilde{E}_i(b_1\otimes b_2)&=\begin{cases}
b_1\otimes \widetilde{E}_ib_2, &\text{if
$\epsilon_i(b_1)\leq\varphi_i(b_2)$;}\\
\widetilde{E}_ib_1\otimes b_2, &\text{if
$\epsilon_i(b_1)>\varphi_i(b_2)$;}\\
\end{cases}\\
\widetilde{F}_i(b_1\otimes b_2)&=\begin{cases}
b_1\otimes \widetilde{F}_ib_2, &\text{if
$\epsilon_i(b_1)<\varphi_i(b_2)$;}\\
\widetilde{F}_ib_1\otimes b_2, &\text{if
$\epsilon_i(b_1)\geq\varphi_i(b_2)$,}\\
\end{cases}
\end{aligned}
$$
where $b_1\in B_1, b_2\in B_2$, $B_1, B_2$ are the crystal bases of
the ${\mmU}_{\mmQ(q)}(\mathfrak{gl}_{m})$-modules $M', M''$
respectively, and $$ \epsilon_i(b)=\max\{n\geq
0|\widetilde{E}_i^nb\neq 0\},\,\,\, \varphi_i(b)=\max\{n\geq
0|\widetilde{F}_i^nb\neq 0\}.
$$
Let $1\leq i\leq m-1$ and $1\leq j\leq m$. Recall that $$
\widetilde{E}_iv_j=\begin{cases} v_i, &\text{if $j=i+1$;}\\
0, &\text{otherwise,} \end{cases}\quad\,
\widetilde{F}_iv_j=\begin{cases} v_{i+1}, &\text{if $j=i$;}\\
0, &\text{otherwise.} \end{cases}
$$
To prove the lemma, it suffices to show that for each $1\leq i\leq m$, $$
\widetilde{E}_i\Bigl(\underbrace{v_{\mu_1}\otimes\cdots\otimes v_2\otimes v_1}_{\text{$\mu_1$ terms}}\otimes
\underbrace{v_{\mu_{2}}\otimes\cdots\otimes v_2\otimes v_1}_{\text{$\mu_{2}$ terms}}\otimes\cdots\otimes
\underbrace{v_{\mu_l}\otimes\cdots\otimes v_2\otimes v_1}_{\text{$\mu_l$ terms}}\Bigr)=0.
$$
But this follows from a direct verification (using the above tensor product rule) and an inductive argument.
\end{proof}

\begin{lemma} \label{keylm00} Let $\mu=(\mu_1,\cdots,\mu_m)$ be a partition of $n$ with $\ell(\mu)\leq m$. Suppose that $b_0\in\hat{B}$ is
a weight vector of weight $\mu$ with respect to the action of $\mathbb{U}_{\mmQ(q)}(\mathfrak{gl}_{m})$. Then

1) $b_0\in B$ and $b_0$ is also a weight vector of weight $\mu$ with respect to the action of $\mathbb{U}_{\mmQ(q)}(\mathfrak{sp}_{2m})$; and

2) if $b_0\in\hat{B}[\mu]^{hi}$ then $b_0\in B[\mu]^{hi}$; and

3) $\hat{M}[>\!\mu]_{\mmA}\subseteq {M}[>\!\mu]_{\mmA}$.
\end{lemma}

\begin{proof} We write $b_0=v_{p_1}\diamond\cdots\diamond v_{p_n}$, where $1\leq p_1,\cdots,p_n\leq m$.
By \cite[Theorem 27.3.2]{Lu}, there is a bar involution $\psi_C$ on $V_{\mmQ(q)}^{\otimes n}$, and for each $1\leq i_1,\cdots,i_n\leq 2m$, there is a unique
canonical basis element $v_{i_1}\diamond_{C}\cdots\diamond_{C} v_{i_n}\in V_{\mmA}^{\otimes n}$ satisfying:
\begin{enumerate}
\item $\psi_{C}(v_{i_1}\diamond_C\cdots\diamond_C v_{i_n})=v_{i_1}\diamond_C\cdots\diamond_C v_{i_n}$, and
\item $v_{i_1}\diamond_C\cdots\diamond_C v_{i_n}$ is equal to $v_{i_1}\otimes\cdots
\otimes v_{i_n}$ plus a linear combination of elements $v_{j_1}\otimes\cdots\otimes v_{j_n}$ with $(v_{j_1},\cdots,v_{j_n})<^{C}(v_{i_1},\cdots,v_{i_n})$, with coefficients in $q^{-1}\mmZ[q^{-1}]$.
\end{enumerate}
Here we use the subscript and superscript ``C" to emphasize that these objects depend on the type $C_{m}$ root system.

We claim that $\psi_C(v_{\mui})=\psi_A(v_{\mui})$ for any
$v_{\mui}\in\hat{V}_{\mmA}^{\otimes n}$. We use induction on $n$.
Suppose that our claim is true when $n$ is replaced by $n-1$. Let
$\psi_A^{(n-1)}$, $\psi_C^{(n-1)}$ be the bar involutions defined on
$\hat{V}_{\mmA}^{\otimes n-1}$, $V_{\mmA}^{\otimes n-1}$
respectively. Let $\Theta$ be the quasi-$\mathcal{R}$-matrix of
${\mmU}_{\mmQ(q)}(\mathfrak{sp}_{2m})$ introduced in \cite[4.1]{Lu}.
We shall mainly follow the notations in \cite[4.1]{Lu} without any
further explanation, except that we use $q$ to denote the quantum
parameter instead of $v$ in Lusztig's book. We define
$$\begin{aligned} &\mathbb{N}:=\mmZ^{\geq
0},\,\,I_0:=\bigl\{\eps_1-\eps_2,\cdots,\eps_{m-1}-\eps_m\bigr\},\,\,I:=\bigl\{\eps_1-\eps_2,
\cdots,\eps_{m-1}-\eps_m,2\eps_m\bigr\},\\
&\mathbf{U}:={\mmU}_{\mmQ(q)}(\mathfrak{sp}_{2m}),\,\, \mathbf{U}^{+}:=\langle e_1,e_2,\cdots,e_m\rangle,\,\,
 \mathbf{U}^{-}:=\langle f_1,f_2,\cdots,f_m\rangle .
\end{aligned}
$$
In a similar way, one can define the subalgebras
${\mmU}_{\mmQ(q)}(\mathfrak{gl}_{m})^{+}$,
${\mmU}_{\mmQ(q)}(\mathfrak{gl}_{m})^{-}$ of
${\mmU}_{\mmQ(q)}(\mathfrak{gl}_{m})$ using the Chevalley generators
of ${\mmU}_{\mmQ(q)}(\mathfrak{gl}_{m})$. It is well-known that both
$\mathbf{U}^+$ and ${\mmU}_{\mmQ(q)}(\mathfrak{gl}_{m})^{+}$ can be
presented by their Chevalley generators and quantum Serre relations
(\cite[4.21]{Ja1}). Using the PBW bases (\cite[8.24]{Ja1}) for
$\mathbf{U}^+$ and ${\mmU}_{\mmQ(q)}(\mathfrak{gl}_{m})^{+}$, it is
easy to see that ${\mmU}_{\mmQ(q)}(\mathfrak{gl}_{m})^{+}$ can be
identified with the subalgebra of $\bf{U}^+$ generated by
$e_1,\cdots,e_{m-1}$. A similar statement holds for
${\mmU}_{\mmQ(q)}(\mathfrak{gl}_{m})^{-}$ and $\mathbf{U}^-$. Let
$\mathbf{f}$ be defined as in \cite[1.2.5]{Lu} and $\{b\}$ be a
$\mmQ(q)$-basis of $\mathbf{f}$ such that
$B_{\nu}:=\{b\}\cap\mathbf{f}_{\nu}$ is a basis of
$\mathbf{f}_{\nu}$ for any $\nu\in\mathbb{N}I$. Let $\{b^*\}$ be the
basis of $\mathbf{f}$ dual to $\{b\}$ under the bilinear form $(\,
,\,)$ introduced in \cite[1.2.12]{Lu}. Then
$\Theta=\sum_{\nu\in\mathbb{N}I}\Theta_{\nu}$, where $$
\Theta_{\nu}=(-1)^{{\rm tr}\nu}q_{\nu}\sum_{b\in
B_{\nu}}b^{-}\otimes
(b^{*})^{+}\in\mathbf{U}_{\nu}^{-}\otimes\mathbf{U}_{\nu}^{+}.
$$
Applying \cite[1.2.3 (d)]{Lu}, we see that $\{b^*\}$ is homogeneous too. Moreover, $\cup_{\substack{\nu\in\mathbb{N}I_0\\ b\in B_{\nu}}}b^{-}$ and $\cup_{\substack{\nu\in\mathbb{N}I_0\\ b\in B_{\nu}}}(b^{*})^{+}$ are bases of ${\mmU}_{\mmQ(q)}(\mathfrak{gl}_{m})^{-}$ and ${\mmU}_{\mmQ(q)}(\mathfrak{gl}_{m})^{+}$ respectively.

Recall that $v_1,\cdots,v_m$ is a part of the canonical bases of the  ${\mmU}_{\mmQ(q)}(\mathfrak{sp}_{2m})$-module $V_{\mmQ(q)}$. We refer the reader to \cite[Section 6]{Hu2} for the description of the whole canonical bases of $V_{\mmQ(q)}$ and the action of ${\mmU}_{\mmQ(q)}(\mathfrak{sp}_{2m})$ on $V_{\mmQ(q)}$. The important facts that we need here are $e_m\hat{V}_{\mmQ(q)}=0$ and $e_j\hat{V}_{\mmQ(q)}\subseteq\hat{V}_{\mmQ(q)}$ for any $1\leq j\leq m-1$.
For any $\nu\in\mathbb{N}I\setminus\mathbb{N}I_0$ and any $b\in B_{\nu}$, we see that $e_m$ must appear in every monomial occurring in $(b^*)^+$. In particular, this implies that $(b^{*})^{+}v_{i_n}=0$ since $1\leq i_n\leq m$.

Therefore, we have that $$\begin{aligned}
\psi_{C}\bigl(v_{\mui}\bigr)&=\Theta\Bigl(\psi_C^{(n-1)}\bigl(v_{i_1}\otimes\cdots\otimes v_{i_{n-1}}\bigr)\otimes\psi_C^{(1)}(v_{i_n})\Bigr)\\
&=\Theta\Bigl(\psi_A^{(n-1)}\bigl(v_{i_1}\otimes\cdots\otimes v_{i_{n-1}}\bigr)\otimes v_{i_n}\Bigr)\\
&=\sum_{\nu\in\mathbb{N}I}(-1)^{{\rm tr}\nu}q_{\nu}\sum_{b\in B_{\nu}}\biggl(b^{-}\psi_A^{(n-1)}\bigl(v_{i_1}\otimes\cdots\otimes v_{i_{n-1}}\bigr)\otimes (b^{*})^{+}v_{i_n}\biggr)\\
&=\sum_{\nu\in\mathbb{N}I_0}(-1)^{{\rm tr}\nu}q_{\nu}\sum_{b\in B_{\nu}}\biggl(b^{-}\psi_A^{(n-1)}\bigl(v_{i_1}\otimes\cdots\otimes v_{i_{n-1}}\bigr)\otimes (b^{*})^{+}v_{i_n}\biggr)\\
&=\psi_A\bigl(v_{\mui}\bigr),
\end{aligned}
$$
as required. This proves our claim. As a result, we deduce that $\psi_C(b_0)=\psi_A(b_0)=b_0$.

Now we have that \begin{equation} \label{bar1}
                  \psi_C\bigl(b_0-v_{p_1}\diamond_C\cdots\diamond_C v_{p_n}\bigr)=b_0-v_{p_1}\diamond_C\cdots\diamond_C v_{p_n}.
                 \end{equation}
On the other hand, by construction, we have that $$
b_0-v_{p_1}\diamond_C\cdots\diamond_C v_{p_n}\in\sum_{\muj\in
I(2m,n)}q^{-1}\mmZ[q^{-1}]v_{j_1}\otimes\cdots\otimes v_{j_n}.
$$
Suppose that $b_0\neq v_{p_1}\diamond_C\cdots\diamond_C v_{p_n}$.
Then we can find an $\mui\in I(2m,n)$, such that
$$
b_0-v_{p_1}\diamond_C\cdots\diamond_C v_{p_n}\in
q^{-1}f(q^{-1})v_{i_1}\diamond_C\cdots\diamond_C
v_{i_n}+\sum_{\substack{\muj\in I(2m,n)\\
\muj\not\ge^{C}\mui}}\mmZ[q,q^{-1}]v_{j_1}\diamond_C\cdots\diamond_C
v_{j_n},
$$
where $0\neq f(q^{-1})\in\mmZ[q^{-1}]$. Applying $\psi_C$ to both sides of the above equality, we get a contradiction to
(\ref{bar1}). This proves that $b_0=v_{p_1}\diamond_C\cdots\diamond_C v_{p_n}\in B$ and the first statement of the lemma follows. The second statement follows from the definition of $B[\mu]^{hi}$ (see \cite[27.2.3]{Lu}), the facts that $\widetilde{E}_mv_i=0$ for any $1\leq i\leq m$ and
the action of the Kashiwara operators $\widetilde{E}_1,\cdots,\widetilde{E}_{m-1}$ are compatible with the natural embedding $\hat{V}_{\mmQ(q)}^{\otimes n}\hookrightarrow{V}_{\mmQ(q)}^{\otimes n}$.

Finally, note that for any $b\in B[\mu]$, $b\in B[\mu]^{hi}$ if and
only if the canonical image of $b$ in
$M[\geq\!\mu]_{\mmA}/M[>\!\mu]_{\mmA}$ is mapped to a highest weight
vector of some direct summand $\Delta(\lam)_{\mmA}$ in the right-hand
side of (\ref{filtra3}).  It follows that $$
{M}[>\!\mu]_{\mmA}=\sum_{\substack{b\in B[\nu]^{hi}\\
\mu<\nu\in\pi_0}}\mathbb{U}_{\mmA}(\mathfrak{sp}_{2m})b,\quad
\hat{M}[>\!\mu]_{\mmA}=\sum_{\substack{b\in\hat{B}[\nu]^{hi}\\
\mu<\nu\in\hat\pi_0}}\mathbb{U}_{\mmA}(\mathfrak{gl}_{m})b.
$$
Now the third statement follows from 2) and the fact for any $v\in\hat{V}_{\mmA}^{\otimes n}$, $\mathbb{U}_{\mmA}(\mathfrak{gl}_{m})v=\mathbb{U}_{\mmA}(\mathfrak{sl}_{m})v$ and every Chevalley generator of $\mathbb{U}_{\mmQ(q)}(\mathfrak{sp}_{2m})$ acts on $v$ in the same way as the corresponding Chevalley generator of $\mathbb{U}_{\mmQ(q)}(\mathfrak{sl}_{m})$.
\end{proof}

Now we are in position to prove the key lemma in this section. The
most difficult part is the second statement of the following lemma.
The main idea of its proof is to show that there exists a canonical
basis element $b\in B[\lam]^{hi}$ of the
${\mmU}_{\mmQ(q)}(\mathfrak{sp}_{2m})$-module $V_{\mmQ(q)}^{\otimes n}$ such that
$b\!\downarrow_{q=1}$ always appears with coefficient $1$ in the linear expansion of $z_{0,\lam}$ into the specialization
at $q=1$ of the canonical basis elements of $V_{\mmQ(q)}^{\otimes n}$. The
element $b$ will be identified with a canonical basis element in the
${\mmU}_{\mmQ(q)}(\mathfrak{gl}_{m})$-module $\hat{V}_{\mmQ(q)}^{\otimes n}$ and
hence eventually identified with a parabolic Kazhdan-Lusztig basis
of certain permutation module over a type $A$ Hecke algebra.

\begin{lemma} \label{lm44} Let $K$ be an algebraically closed field. Suppose that $\lam$ is  a
partition of $n$ with $\ell(\lam)\leq m$. Then
\begin{enumerate}
\item $z_{0,\lam}$ is a non-zero maximal vector of weight $\lam$ in $V^{\otimes
n}$;
\item $(KG)z_{0,\lam}\cong\Delta(\lam)$, and $V^{\otimes n}/(KG)z_{0,\lam}$ has a Weyl filtration;
\item the dimension of $z_{0,\lam}\mbb_n$ is independent of the
field $K$.
\end{enumerate}
\end{lemma}

\begin{proof} It is well-known that $v_{\lam}w_{\lam}x_{\lam'}\neq 0$ (cf. \cite{Gr}).
In particular, $z_{0,\lam}=v_{\lam}w_{\lam}x_{\lam'}\neq 0$. Note that $z_{0,\lam}=Z_{0,\lam}\!\downarrow_{q=1}$. Thus the statement 1) follows from Lemma \ref{lmhx}. Since $z_{0,\lam}\in \hat{V}^{\otimes n}$, it follows easily that
$z_{0,\lam}\mbb_n=z_{0,\lam}K\bBS_n$. It is also well-known (cf. \cite{Gr}) that
$z_{0,\lam}KS_n\cong y_{\lam}w_{\lam}x_{\lam'}K\bBS_n$ as a right $K\bBS_n$-module. In particular, $$
\dim z_{0,\lam}\mbb_n=\#\bigl\{\text{standard $\lam'$-tableaux}\bigr\}, $$ which is independent of the field $K$. This
proves the statement 3). It remains to prove the statement 2). We divide the proof into
three steps: \smallskip

{\it Step 1.} We claim that there exists a canonical basis element
${b}\in \hat{B}[\lam]^{hi}$, such that $$ z_{0,\lam}\equiv
{b}\!\!\downarrow_{q=1}+\sum_{b\neq b'\in
\hat{B}[\lam]^{hi}}{c}'_{b'}b'\!\!\downarrow_{q=1}\!\!\!\pmod{\hat{M}[>\!\lam]_{\mmZ}},
$$
where ${c}'_{b'}\in\mmZ$ for each $b'$.

For each composition $\mu$ of $n$, let
$\bigl(\hat{V}_{\mmA}^{\otimes n}\bigr)_{\mu}$ be the corresponding
weight subspace of $\hat{V}_{\mmA}^{\otimes n}$ (with respect to the
action of $\mathbb{U}_{\mmA}(\mathfrak{gl}_{m})$). It is well known
that there is a right $\mHH_{\mmA}(\bBS_n)$-module isomorphism:
$$\begin{aligned}
\phi_{\mu}:\,\,\hat{Y}_{\mu}\mHH_{\mmA}(\bBS_n)&\cong v_{\mu}\mHH_{\mmA}(\bBS_n)=\bigl(\hat{V}_{\mmA}^{\otimes n}\bigr)_{\mu}\cong \Ind_{\mHH_q(\bBS_{{\lam}})}^{\mHH_q(\bBS_n)}\rho_{\mu}\\
\hat{Y}_{\mu}h & \rightarrow v_{\mu}h ,\quad\forall\,h\in \mHH_{\mmA}(\bBS_n),
\end{aligned}
$$
where $\rho_{\mu}$ is the one dimensional representation of $\mHH_q(\bBS_{\mu})$ which is defined on generators by $\hat{T}_i\mapsto -q^{-1}$ for each $s_i\in\bBS_{\mu}$.

Let $\phi$ be the anti-linear involution on $\mHH_q(\bBS_n)$ which is defined on generators by
$\phi(\hat{T}_i)=\hat{T}_i^{-1}, \phi(q)=q^{-1}$ for each $1\leq i\leq n-1$. By the main result of \cite{FKK}, $\phi$ naturally induces an anti-linear involution on $\hat{Y}_{\lam}\mHH_q(\bBS_n)$ so that one can define the parabolic Kazhdan-Lusztig bases $\bigl\{C_d\bigr\}_{d\in\mathcal{D}_{\lam}}$ of $\hat{Y}_{\lam}\mHH_q(\bBS_n)$. Recall that $\hat{Y}_{\lam}\mHH_q(\bBS_n)$ also has a standard basis $\bigl\{\hat{Y}_{\lam}\hat{T}_d\bigr\}_{d\in\mathcal{D}_{\lam}}$.
For each $d\in\mathcal{D}_{\lam}$, we have that \begin{equation}\label{can1}
C_d\equiv\hat{Y}_{\lam}\hat{T}_d\!\!\pmod{\sum_{d>d'\in\mathcal{D}_{\lam}}q^{-1}\mmZ[q^{-1}]\hat{Y}_{\lam}\hat{T}_{d'}},
\end{equation}
where ``$>$" is the usual Bruhat order defined on the symmetric
group $\bBS_n$. We identify $\hat{Y}_{\lam}\mHH_{\mmA}(\bBS_n)$ with
$\bigl(\hat{V}_{\mmA}^{\otimes n}\bigr)_{\lam}$ via the isomorphism
$\phi_{\lam}$. It follows from \cite[Theorem 2.5]{FKK} that the
bases $\bigl\{C_d\bigr\}_{d\in\mathcal{D}_{\lam}}$ coincide with the
canonical bases of $\hat{V}_{\mmQ(q)}^{\otimes n}$ which are of
weight $\lam$. Here one should understand the notations
$q,\hat{T}_i$ in this paper as the notations $v, -vT_i$ in
\cite{FKK}. Let $\mui_{\lam}\in I(m,n)$ such that
$v_{\mui_{\lam}}=v_{\lam}$. For each $d\in\mathcal{D}_{\lam}$,
$\phi_{\lam}\bigl(\hat{Y}_{\lam}\hat{T}_d\bigr)=v_{\mui_{\lam}d}$.
Note that if $\mui_{\lam}d=(j_1,\cdots,j_n)$, then
$\phi_{\lam}(C_d)=v_{j_1}\diamond v_{j_2}\diamond\cdots\diamond
v_{j_n}$.

We want to use the isomorphism $\phi_{\lam}$ and (\ref{can1}) to
express $z_{0,\lam}$ into a linear combination of the
specializations at $q=1$ of some canonical basis elements. To this
end, we first express $z_{0,\lam}$ into a linear combination of
standard basis elements. Let $w_{0,\lam'}$ be the longest element in
the Young subgroup $\bBS_{\lam'}$ of $\bBS_n$. Note that
$\mft^{\lam}w_{\lam}w_{0,\lam'}$ is row-standard, which implies that
$w_{\lam}w_{0,\lam'}\in\mathcal{D}_{\lam}$. Since
$w_{\lam}\in\mathcal{D}_{\lam}\bigcap\mathcal{D}_{\lam'}^{-1}$, we
have that $\ell(w_{\lam}w_{0,\lam'})=
\ell(w_{\lam})+\ell(w_{0,\lam'})$. Note also that
$w_{\lam}^{-1}\bBS_{\lam}w_{\lam}\bigcap\bBS_{\lam'}=\{1\}$. It
follows that $$
y_{\lam}w_{\lam}x_{\lam'}=y_{\lam}w_{\lam}w_{0,\lam'}+\sum_{\substack{d\in\mathcal{D}_{\lam}\\
d\neq w_{\lam}w_{0,\lam'}}}a_{d}y_{\lam}d
$$
where $a_d\in\mmZ$ for each $d$ and $y_{\lam}:=\sum_{w\in\bBS_{\lam}}(-1)^{\ell(w)}w$. Since $w_{\lam}w_{0,\lam'}>w_{\lam}w$
for any $w_{0,\lam'}\neq w\in\bBS_{\lam'}$, we deduce that $a_d\neq 0$ only if $d<w_{\lam}w_{0,\lam'}$. That is,
\begin{equation}\label{can2}
y_{\lam}w_{\lam}x_{\lam'}=y_{\lam}w_{\lam}w_{0,\lam'}+\sum_{\substack{d\in\mathcal{D}_{\lam}\\
d<w_{\lam}w_{0,\lam'}}}a_{d}y_{\lam}d. \end{equation} Recall that in
this paper $\bBS_n$ acts on $V^{\otimes n}$ by sign permutation
action. Thus, $v_{\lam}\mmZ\bBS_n\cong y_{\lam}\mmZ\bBS_n$,
$v_{\lam}h\mapsto y_{\lam}h,\,\forall\,h\in\mmZ\bBS_n$. Note that
$y_{\lam}=(\hat{Y}_{\lam})\!\downarrow_{q=1}$, $d=(\hat{T}_d)\!\downarrow_{q=1}$
for each $d\in\mathcal{D}_{\lam}$. It follows from (\ref{can1}),
(\ref{can2})  and the isomorphism $\phi_{\lam}$ that
\begin{equation}\label{expan1}
v_{\lam}w_{\lam}x_{\lam'}=b\!\downarrow_{q=1}+\sum_{b\neq b'\in
\hat{B}}c_{b'}b',
\end{equation}
where $b:=\phi_{\lam}(C_{w_{\lam}w_{0,\lam'}})\in\hat{B}$. We set
$\mu:=\lam'$ and write $\mu=(\mu_1,\cdots,\mu_l)$, where $\mu_l\neq
0$ and $\mu_1=\ell(\lam)\leq m$. Then it is easy to check that $$
\mui_{\lam}(w_{\lam}w_{0,\lam'})=\bigl(\underbrace{\mu_1,\mu_1-1,\cdots,1}_{\text{$\mu_1$
terms}},\underbrace{\mu_2, \mu_2-1,\cdots,1}_{\text{$\mu_2$
terms}},\cdots,\underbrace{\mu_l,\mu_l-1,\cdots,1}_{\text{$\mu_l$
terms}}\bigr),$$ which implies that $b=\widehat{b}_{\mu}\in
\hat{B}^{hi}$ by Lemma \ref{lmkey41}. Hence $b\in
\hat{B}[\lam]^{hi}$ as $b$ is a weight vector of weight $\lam$. Note
that $c'_{b'}\neq 0$ only if $b'$ is a weight vector of weight
$\lam$. On the other hand, if $b'\in\hat{B}[\mu]$ and $b'$ is a
weight vector of weight $\lam$ then we must have that $\lam$ is a
weight of the Weyl module of ${\mmU}_{\mmQ(q)}(\mathfrak{gl}_{m})$
associated to $\mu$, which implies that $\mu\geq_{A}\lam$.
Therefore, we can rewrite (\ref{expan1}) as  $$ z_{0,\lam}\equiv
{b}\!\!\downarrow_{q=1}+\sum_{b\neq b'\in
\hat{B}[\lam]^{hi}}{c}'_{b'}b'\!\!\downarrow_{q=1}\!\!\!\pmod{\hat{M}[>\!\lam]_{\mmZ}},
$$
where $b\in \hat{B}[\lam]^{hi}$, ${c}'_{b'}\in\mmZ$ for each $b'$.
This proves our claim.
\smallskip

{\it Step 2.} We claim that \begin{equation}\label{filt6}
z_{0,\lam}\equiv b\!\!\downarrow_{q=1}+\sum_{b\neq b'\in
B[\lam]^{hi}}c_{b'}b'\!\!\downarrow_{q=1}\!\!\!\pmod{M[>\!\lam]_{\mmZ}},
\end{equation}
where $b\in B[\lam]^{hi}$ and $c_{b'}\in\mmZ$ for each $b'$.

Indeed, this follows directly from the result we obtained in Step 1 and Lemma \ref{keylm00}.
\smallskip

{\it Step 3.} Let $\mathbb{U}_{\mmZ}(\mathfrak{sp}_{2m})$ be the
Kostant $\mmZ$-form in $\mathbb{U}_{\mmQ}(\mathfrak{sp}_{2m})$. It is
well-known that (cf. \cite[(6.7)(c),(6.7)(d)]{Lu1}) $$
\mathbb{U}_{\mmZ}(\mathfrak{sp}_{2m})\cong\Bigl(\mathbb{U}_{\mmA}(\mathfrak{sp}_{2m})\otimes_{\mmA}\mmZ\Bigr)/\bigl(\langle
K_1-1,\cdots,K_m-1\rangle\otimes_{\mmA}\mmZ\bigr).
$$
We claim that
$\mathbb{U}_{\mmZ}(\mathfrak{sp}_{2m})z_{0,\lam}\cong\Delta_{\mmZ}(\lam)$,
and $V_{\mmZ}^{\otimes
n}/\mathbb{U}_{\mmZ}(\mathfrak{sp}_{2m})z_{0,\lam}$ has a Weyl
filtration.

Recall the equality (\ref{filt6}) we obtained in Step 2. Let $b\in
B[\lam]^{hi}$ be the canonical basis element we obtained in Step 2.
Using Lemma 3.7 3), we can get a sequence $\{\mu_{(i)}\}_{i=1}^{k}$
of dominant weights in $\pi_0$ and construct a Weyl filtration of
$M_{\mmA}$ as follows:
\begin{equation}\label{filtra4} 0=M^{\mmA}_{(0)}\subset
M^{\mmA}_{(1)}\subset\cdots\subset M^{\mmA}_{(k-1)}\subset
M^{\mmA}_{(k)}=M_{\mmA}=V_{\mmA}^{\otimes n},
\end{equation}
such that \begin{enumerate}
\item for each integer $1\leq i\leq k$, $M^{\mmA}_{(i)}$ is spanned by the canonical basis elements it contains;
\item for each integer $1\leq i\leq k$, there is a ${\mmU}_{\mmA}(\mathfrak{sp}_{2m})$-module
isomorphism:
$M^{\mmA}_{(i)}/M^{\mmA}_{(i-1)}\cong\Delta(\mu_{(i)})_{\mmA}$, such that
if $b'\in M^{\mmA}_{(i)}\cap B$ then the canonical image of $b'$ in
$M^{\mmA}_{(i)}/M^{\mmA}_{(i-1)}$ is mapped either to $0$ or to a
canonical basis element of $\Delta(\mu_{(i)})_{\mmA}$;
\item $\mu_{(i)}>\mu_{(j)}$ only if $i<j$;
\item $B[\lam]^{hi}\setminus\{b\}\subseteq M^{\mmA}_{(i_0-1)}$, where $1\leq i_0\leq k$ is the unique integer such that
$b\in M^{\mmA}_{(i_0)}\setminus M^{\mmA}_{(i_0-1)}$.
\end{enumerate} In particular, $\mu_{(i_0)}=\lam$ as $b\in B[\lam]^{hi}$. Thus
we can write $B=\{b_i\}_{i=1}^N$, where $N=(2m)^n$, such that for
each $1\leq i\leq k$, the elements in the subset
$\bigl\{b_j\bigm|1\leq
j\leq\sum_{s=1}^i\dim\Delta(\mu_{(s)})\bigr\}$ form an $\mmA$-basis of
$M^{\mmA}_{(i)}$. Moreover,
\begin{equation}\label{generator} M^{\mmA}_{(i)}=\sum_{b\in M^{\mmA}_{(i)}\cap B^{hi}}{\mmU}_{\mmA}(\mathfrak{sp}_{2m})b.
\end{equation}
Specializing at $q=1$, we get a
${\mmU}_{\mmZ}(\mathfrak{sp}_{2m})$-submodules filtration of
$V_{\mmZ}^{\otimes n}$. We define a new
${\mmU}_{\mmZ}(\mathfrak{sp}_{2m})$-submodules filtration of
$V_{\mmZ}^{\otimes n}$ as follows: \begin{equation}\label{filtra5}
0=N^{\mmZ}_{(0)}\subset N^{\mmZ}_{(1)}\subset\cdots\subset
N^{\mmZ}_{(k-1)}\subset N^{\mmZ}_{(k)}=M_{\mmZ}=V_{\mmZ}^{\otimes n},
\end{equation}
where $$ N^{\mmZ}_{(i)}:=\begin{cases}
{\mmU}_{\mmZ}(\mathfrak{sp}_{2m})z_{0,\lam}+\sum\limits_{b'\in M^{\mmA}_{(i-1)}\cap
B^{hi}}{\mmU}_{\mmZ}(\mathfrak{sp}_{2m})b'\!\!\downarrow_{q=1},
&\text{if $1\leq i<i_0$;}\\
M^{\mmZ}_{(i)}, &\text{if $i_0\leq i\leq k$.}
\end{cases}
$$
We claim that (\ref{filtra5}) is a Weyl filtration of $M_{\mmZ}$. In
fact, to prove this claim, it suffices to show that for each $1\leq
i\leq i_0$,
$N^{\mmZ}_{(i)}/N^{\mmZ}_{(i-1)}$ is isomorphic to some Weyl module.

Let $F$ be any field which is an $\mmZ$-algebra. Let ${\mmU}_{F}(\mathfrak{sp}_{2m}):={\mmU}_{\mmZ}(\mathfrak{sp}_{2m})\otimes_{\mmZ}F$. Replacing $\mmZ$ by $F$ in (\ref{filtra5}), we get a
${\mmU}_{F}(\mathfrak{sp}_{2m})$-submodules filtration of
$M_{F}:=V_{\mmZ}^{\otimes n}\otimes F$ as follows: \begin{equation}\label{filtra51}
0=N^{F}_{(0)}\subset N^{F}_{(1)}\subset\cdots\subset
N^{F}_{(k-1)}\subset N^{F}_{(k)}=M_{F}.
\end{equation}
Recall that $z_{0,\lam}$ is a maximal vector of weight $\lam$. Using
(\ref{filt6}), (\ref{generator}), together with the third and the
fourth properties of the filtration (\ref{filtra4}), we see that
$$
N^{F}_{(i_0)}={\mmU}_{F}(\mathfrak{sp}_{2m})z_{0,\lam}+\sum\limits_{b'\in M^{\mmA}_{(i_0-1)}\cap
B^{hi}}{\mmU}_{F}(\mathfrak{sp}_{2m})b'\!\!\downarrow_{q=1}.
$$
By the second property of the filtration (\ref{filtra4}), we know that $$
\biggl(\sum\limits_{b'\in M^{\mmA}_{(i)}\cap B^{hi}}{\mmU}_{F}(\mathfrak{sp}_{2m})b'\!\!\downarrow_{q=1}\biggr)/\biggl(\sum\limits_{b'\in M^{\mmA}_{(i-1)}\cap B^{hi}}{\mmU}_{F}(\mathfrak{sp}_{2m})b'\!\!\downarrow_{q=1}\biggr)\cong\Delta(\mu_{(i)})_{F}.
$$
It follows that each $N^{F}_{(i)}/N^{F}_{(i-1)}$ must be a
homomorphic image of $\Delta(\lam_{(i)})_{F}$, where $$
\lam_{(i)}:=\begin{cases} \lam, &\text{if $i=1$;}\\
\mu_{(i-1)}, &\text{$1\leq i\leq i_0$;}\\
\mu_{(i)}, &\text{if $i_0+1\leq i\leq k$.}
\end{cases}
$$ Noting that $$
\dim(M_{F})=(2m)^n=\sum_{i=1}^k\rank\Delta(\mu_{(i)})_{\mmZ}=\sum_{i=1}^k\dim\Delta(\mu_{(i)})_{F}=\sum_{i=1}^k\dim\Delta(\lam_{(i)})_{F},
$$
and comparing the dimensions, we deduce that the natural surjection from
$\Delta(\lam_{(i)})_{F}$ onto
$N^{F}_{(i)}/N^{F}_{(i-1)}$ must be an isomorphism. In particular, $$
\mathbb{U}_{F}(\mathfrak{sp}_{2m})z_{0,\lam}=N^{F}_{(1)}\cong \Delta(\lam_{(1)})_{F}
=\Delta(\lam)_{F}. $$
Since $z_{0,\lam}$ is a maximal vector of weight $\lam$, we have a natural surjection from
$\Delta(\lam)_{\mmZ}$ onto $\mathbb{U}_{\mmZ}(\mathfrak{sp}_{2m})z_{0,\lam}$, which induces a surjection from
$\Delta(\lam)_{\mmZ}\otimes F$ onto $\mathbb{U}_{\mmZ}(\mathfrak{sp}_{2m})z_{0,\lam}\otimes F$. On the other hand,
we also have a natural surjection from $\mathbb{U}_{\mmZ}(\mathfrak{sp}_{2m})z_{0,\lam}\otimes F$ onto
$\mathbb{U}_{F}(\mathfrak{sp}_{2m})z_{0,\lam}\cong\Delta(\lam)_{F}$. It follows again by comparing dimensions that the natural surjection from
$\Delta(\lam)_{\mmZ}\otimes F$ onto $\mathbb{U}_{\mmZ}(\mathfrak{sp}_{2m})z_{0,\lam}\otimes F$ is an isomorphism for any
field $F$. Hence the natural surjection from
$\Delta(\lam)_{\mmZ}$ onto $\mathbb{U}_{\mmZ}(\mathfrak{sp}_{2m})z_{0,\lam}$ must be an isomorphism as well. By a similar
argument, we can show that the natural surjection from $\Delta(\lam_{(i)})_{\mmZ}$ onto $N^{\mmZ}_{(i)}/N^{\mmZ}_{(i-1)}$
must be an isomorphism as well for each $i$. As a result, the $\mathbb{U}_{\mmZ}(\mathfrak{sp}_{2m})$-module $V_{\mmZ}^{\otimes
n}/\mathbb{U}_{\mmZ}(\mathfrak{sp}_{2m})z_{0,\lam}=N^{\mmZ}_{(k)}/N^{\mmZ}_{(1)}$ has a Weyl
filtration. This proves our claim.

Finally, by taking $F=K$ and noting that
$\mathbb{U}_{K}(\mathfrak{sp}_{2m})z_{0,\lam}=(KG)z_{0,\lam}$ (cf.
\cite{Ja}) whenever $K$ is an algebraically closed field, we deduce
that $(KG)z_{0,\lam}\cong\Delta(\lam)$, and $V^{\otimes
n}/(KG)z_{0,\lam}$ has a Weyl filtration. This completes the proof
of the statement 2).
\end{proof}

\begin{corollary} \label{cor45} Let $g$ be an integer with $0\leq g\leq [n/2]$ and $\lam$ a partition of $n-2g$ satisfying
$\ell(\lam)\leq m$. Then there exists an embedding
$\Delta(\lam)\hookrightarrow
V^{\otimes n-2g}$ such that $V^{\otimes n-2g}/\Delta(\lam)$ has a Weyl filtration. In particular, $$
\Ext_{G}^{1}\bigl(V^{\otimes n-2g}/\Delta(\lam),V^{\otimes n}\bigr)=0.
$$
\end{corollary}

\begin{proof} This follows from Lemma \ref{lm44} and Lemma \ref{BP}.
\end{proof}

In the remaining part of this section, we fix an integer $0\leq
g\leq [n/2]$ and a partition $\lam$ of $n-2g$ with $\ell(\lam)\leq
m$. For simplicity, we shall write $z_{g,\lam}$ instead of
$z_{g,\lam}\otimes_{\mmZ}1_{K}$. By Lemma \ref{lmhx} and specializing
$q$ to $1$, we get that $z_{g,\lam}$ is a non-zero maximal vector of
weight $\lam$ with respect to the action of $G$ on
$V^{\otimes n}$. That is,
\begin{equation}\label{mvector}
0\neq z_{g,\lam}\in \bigl(V^{\otimes n}\bigr)^{U}_{\lam}.
\end{equation}
As a consequence, $z_{g,\lam}\mbb_n\subseteq \bigl(V^{\otimes n}\bigr)^{U}_{\lam}$.
On the other hand, we have $$
\bigl(V^{\otimes n}\bigr)^{U}_{\lam}\cong\Hom_{KG}\bigl(\Delta(\lam),V^{\otimes n}\bigr).
$$
Since $V^{\otimes n}$ has a good filtration, it follows from Lemma \ref{BP} that the dimensions of
$\Hom_{KG}\bigl(\Delta_{\lam},V^{\otimes n}\bigr)$ and hence of $\bigl(V^{\otimes n}\bigr)^{U}_{\lam}$ are
independent of
$K$. Therefore, to complete the proof of Lemma \ref{Hu2lem} as well as the second part of Theorem \ref{mainthm3},
it suffices to prove that \begin{equation}\label{equaZ1}
\bigl(V^{\otimes n}\bigr)^{U}_{\lam}=z_{g,\lam}\mbb_n.
\end{equation}

\begin{lemma} \label{lmk40} With the notations as above, we have that $z_{g,\lam}\mbb_n=z_{g,\lam}K\bBS_n$.
\end{lemma}

\begin{proof} For each $1\leq s<t\leq n$ and $\mui\in I(2m,n)$, recall that $$\begin{aligned}
v_{\mui}e_{s,t}&=-\langle v_{i_s}, v_{i_t}\rangle \sum_{k=1}^{2m}v_{i_1}\otimes\cdots\otimes
v_{i_{s-1}}\otimes{v_k}\otimes v_{i_{s+1}}\otimes\cdots\otimes\\
&\qquad\qquad v_{i_{t-1}}\otimes{v_k^{\ast}}\otimes
v_{i_{t+1}}\otimes\cdots\otimes v_{i_n}.\end{aligned}
$$
The lemma follows from the definition of $z_{g,\lam}$ and some direct verification.
\end{proof}
\medskip

In view of the above discussion and Lemma \ref{lmk40}, to prove (\ref{equaZ1}), it suffices to show that
$\bigl(V^{\otimes n}\bigr)^{U}_{\lam}=z_{g,\lam}K\bBS_n$. Using Corollary \ref{cor45}, we have
an embedding $\Delta(\lam)\hookrightarrow
V^{\otimes n-2g}$ such that $V^{\otimes n-2g}/\Delta(\lam)$ has a Weyl filtration.
Therefore, we have the following commutative diagram of
homomorphisms:
\begin{equation}\label{cd3}
\xymatrix{0\ar[r]  & \Hom_{KG}\bigl(V^{\otimes n-2g}, V^{\otimes n}\mbb_n^{(g)}\bigr)
\ar[r]^{\,\quad\sim}\ar[d]_{\beta} & \Hom_{KG}\bigl(V^{\otimes n-2g}, V^{\otimes n}\bigr) \ar@{>>}[d]\\
0 \ar[r] & \Hom_{KG}\bigl(\Delta(\lam), V^{\otimes n}\mbb_n^{(g)}\bigr)
\ar@{^{(}->}[r] & \Hom_{KG}\bigl(\Delta(\lam), V^{\otimes n}\bigr),}
\end{equation}
where by Lemma \ref{keylem1} the top horizontal map is an
isomorphism and the fact that $\Ext_{KG}^{1}\bigl(V^{\otimes
n-2g}/\Delta(\lam),V^{\otimes n}\bigr)=0$ forces that the right
vertical map is a surjection. Since the bottom horizontal map is an
injection, it follows that the left vertical map $\beta$ must be a
surjection and the bottom horizontal map must be an isomorphism. Now
applying Lemma \ref{keylem1},  we get that
$\Hom_{KG}\bigl(\Delta(\lam), V^{\otimes n}\mbb_n^{(g)}\bigr)$ is
spanned by $\beta\bigl(\sigma\tau_g\bigr)$ for $\sigma\in\bBS_n$,
where $\tau_g$ is defined in (\ref{tauf}). Therefore, we can deduce
that the subspace of maximal vectors of weight $\lam$ in $V^{\otimes
n}$ is spanned by all
$\beta\bigl(\sigma\tau_g\bigr)(z_{0,\lam})=z_{g,\lam}\sigma$, where
$\sigma\in\bBS_n$. Hence $\bigl(V^{\otimes n}\bigr)^{U}_{\lam}=
z_{g,\lam}K\bBS_n$. This completes the proof of Lemma \ref{Hu2lem} as
well as the second part of Theorem \ref{mainthm3}.

\section{Proof of Theorem \ref{mainthm2}}

The purpose of this section is to give a proof of Theorem
\ref{mainthm2}.\smallskip

\begin{lemma} \label{hom0} For any integer $1\leq f\leq [n/2]$, we have that $$
\Hom_{KG}\bigl(V^{\otimes n}\mbb_n^{(f)}, V^{\otimes n}/V^{\otimes
n}\mbb_n^{(f)}\bigr)=0.
$$
In particular, the canonical embedding $$\iota_1:
\End_{KG}\bigl(V^{\otimes n}/V^{\otimes
n}\mbb_n^{(f)}\bigr)\hookrightarrow\Hom_{KG}\bigl(V^{\otimes
n},V^{\otimes n}/V^{\otimes n}\mbb_n^{(f)}\bigr)$$ is actually an
isomorphism.
\end{lemma}

\begin{proof} Suppose that $\Hom_{KG}\bigl(V^{\otimes n}\mbb_n^{(f)}, V^{\otimes
n}/V^{\otimes n}\mbb_n^{(f)}\bigr)\neq 0$. By the proof of Lemma
\ref{Jakey2}, $V^{\otimes n}\mbb_n^{(f)}$ is a sum of some submodules
$M_i$ such that $M_i\cong V^{\otimes n-2f}$ for each $i$. It follows
that $$ \Hom_{KG}\bigl(V^{\otimes n-2f}, V^{\otimes n}/V^{\otimes
n}\mbb_n^{(f)}\bigr)\neq 0.
$$
Since $V^{\otimes n-2f}$ is a tilting module, $V^{\otimes n-2f}$ has
a Weyl filtration. For each $\mu\in X^{+}$, $\bigl(V^{\otimes
n-2f}:\Delta(\mu)\bigr)\neq 0$ only if $\mu\vdash (n-2f-2t)$ for
some integer $0\leq t\leq [(n-2f)/2]$.

By Theorem \ref{mainthm1}, we know that $V^{\otimes n}/V^{\otimes
n}\mbb_n^{(f)}$ has a good filtration. For each $\lam\in X^{+}$,
$\bigl(V^{\otimes n}/V^{\otimes n}\mbb_n^{(f)}:\nabla(\lam)\bigr)\neq
0$ only if $\lam\vdash (n-2s)$ for some integer $0\leq s<f$.

By Lemma \ref{BP}, we get that $$\begin{aligned}
&\quad\,\dim\Hom_{KG}\bigl(V^{\otimes n-2f}, V^{\otimes
n}/V^{\otimes n}\mbb_n^{(f)}\bigr)\\
&=\sum_{\lam\in X^{+}}\bigl(V^{\otimes n-2f}:\Delta(\lam)\bigr)
\bigl(V^{\otimes n}/V^{\otimes n}\mbb_n^{(f)}:\nabla(\lam)\bigr)\\
&=\sum_{\substack{\lam\vdash n-2f-2t,\, \lam\vdash n-2s\\
0\leq t\leq [(n-2f)/2]\\
0\leq s<f}}\bigl(V^{\otimes n-2f}:\Delta(\lam)\bigr)
\bigl(V^{\otimes n}/V^{\otimes n}\mbb_n^{(f)}:\nabla(\lam)\bigr)\\
&=0,
\end{aligned}
$$
which is a contradiction.
\end{proof}

\begin{lemma} \label{lm42} 1) The canonical map $$
\theta_1: \End_{KG}\bigl(V^{\otimes
n}\bigr)\rightarrow\Hom_{KG}\bigl(V^{\otimes n},V^{\otimes
n}/V^{\otimes n}\mbb_n^{(f)}\bigr)
$$
is surjective;

2) the dimension of $\End_{KG}\bigl(V^{\otimes n}/V^{\otimes
n}\mbb_n^{(f)}\bigr)$ is independent of $K$.
\end{lemma}

\begin{proof} We have the following exact sequence of maps:
$$\begin{aligned}
&0\rightarrow\Hom_{KG}\bigl(V^{\otimes n},V^{\otimes
n}\mbb_n^{(f)}\bigr)\rightarrow\End_{KG}\bigl(V^{\otimes
n}\bigr)\overset{\theta_1}{\rightarrow}\\
&\qquad \Hom_{KG}\bigl(V^{\otimes n},V^{\otimes n}/V^{\otimes
n}\mbb_n^{(f)}\bigr)\rightarrow\Ext^{1}_{G}\bigl(V^{\otimes
n},V^{\otimes n}\mbb_n^{(f)}\bigr).
\end{aligned}$$
Since $V^{\otimes n}$ has a Weyl filtration, and by Theorem
\ref{mainthm2}, $V^{\otimes n}\mbb_n^{(f)}$ has a good filtration, it
follows that $\Ext^{1}_{G}\bigl(V^{\otimes n},V^{\otimes
n}\mbb_n^{(f)}\bigr)=0$. This implies that the canonical map $$
\theta_1: \End_{KG}\bigl(V^{\otimes
n}\bigr)\rightarrow\Hom_{KG}\bigl(V^{\otimes n},V^{\otimes
n}/V^{\otimes n}\mbb_n^{(f)}\bigr)
$$
is surjective. This proves 1).

The above exact sequence implies that $$ \begin{aligned}&\quad \dim
\Hom_{KG}\bigl(V^{\otimes n},V^{\otimes n}/V^{\otimes
n}\mbb_n^{(f)}\bigr)\\
&=\dim \End_{KG}\bigl(V^{\otimes n}\bigr)-\dim
\Hom_{KG}\bigl(V^{\otimes n},V^{\otimes
n}\mbb_n^{(f)}\bigr).\end{aligned}
$$
Since $V^{\otimes n}$ has a Weyl filtration as well as a good
filtration, and $V^{\otimes n}\mbb_n^{(f)}$ has a good filtration,
and both the character formula of $V^{\otimes n}$ and of $V^{\otimes
n}\mbb_n^{(f)}$ are independent of $K$, it follows that
$$\begin{aligned} &\quad\, \dim \Hom_{KG}\bigl(V^{\otimes n},V^{\otimes
n}/V^{\otimes
n}\mbb_n^{(f)}\bigr)\\
&=\sum_{\lam\in X^{+}}\bigl(V^{\otimes
n}:\Delta(\lam)\bigr)\bigl(V^{\otimes n}:\nabla(\lam)\bigr)
-\sum_{\lam\in X^{+}}\bigl(V^{\otimes
n}:\Delta(\lam)\bigr)\bigl(V^{\otimes
n}\mbb_n^{(f)}:\nabla(\lam)\bigr)\\
&=\sum_{\lam\in X^{+}}\Bigl(\bigl(V_{\mmC}^{\otimes
n}:\Delta_{\mmC}(\lam)\bigr)\bigl(V_{\mmC}^{\otimes
n}:\nabla_{\mmC}(\lam)\bigr)-\\
&\qquad\qquad\qquad\qquad\bigl(V_{\mmC}^{\otimes
n}:\Delta_{\mmC}(\lam)\bigr)\bigl(V_{\mmC}^{\otimes
n}\mbb_n^{(f)}:\nabla_{\mmC}(\lam)\bigr)\Bigr),\end{aligned}
$$
which is independent of $K$. Note that $$
\Hom_{KG}\bigl(V^{\otimes n},V^{\otimes n}/V^{\otimes n}\mbb_n^{(f)}\bigr)=
\End_{KG}\bigl(V^{\otimes n}/V^{\otimes n}\mbb_n^{(f)}\bigr).
$$
So 2) also follows.
\end{proof}
\medskip

\noindent {\bf Proof of Theorem \ref{mainthm2}:} Let $\varphi''_K$
denote the natural $K$-algebra homomorphism:
$\mbb_n\rightarrow\End_{KG}\bigl(V^{\otimes n}/V^{\otimes
n}\mbb_n^{(f)}\bigr)$. Then $$
\varphi_K\bigl(\mbb_n/\mbb_n^{(f)}\bigr)=\varphi''_K(\mbb_n).
$$
In view of Lemma \ref{lm42}, it suffices to show that $\varphi''_K$
is surjective. We consider the following commutative diagram of
maps:
$$
\xymatrix{
 \mbb_n \ar@{>>}[r]^{\varphi_K}\ar[d]_{\varphi''_K} & \End_{KG}\bigl(V^{\otimes n}\bigr)
\ar@{>>}[d]^{\theta_1}\\
 \End_{KG}\bigl(V^{\otimes n}/V^{\otimes
n}\mbb_n^{(f)}\bigr) \ar[r]_{\iota_1\,\,\,\quad}^{\sim\,\,\quad} &
\Hom_{KG}\bigl(V^{\otimes n},V^{\otimes n}/V^{\otimes
n}\mbb_n^{(f)}\bigr),}
$$
By \cite{DDH}, the top horizontal map is surjective. By Lemma
\ref{lm42}, $\theta_1$ is surjective. By Lemma \ref{hom0}, $\iota_1$
is an isomorphism. It follows that $\varphi''_K$ must be surjective,
as required.
\medskip

We remark that Theorem \ref{mainthm2} can be strengthened in the
following sense: the algebraically closed field $K$ can be replaced
by an arbitrary infinite field.

\begin{proposition} \label{mainprop} For any infinite field $K$ (not necessarily algebraically
closed), $\dim_{K}\End_{KSp_{2m}(K)}\bigl(V^{\otimes n}/V^{\otimes
n}\mbb_n^{(f)}\bigr)$ is independent of the infinite field $K$.
Moreover, $$ \varphi_K(\mbb_n/\mbb_n^{(f)})=\End_{
KSp_{2m}(K)}\bigl(V^{\otimes n}/V^{\otimes n}\mbb_n^{(f)}\bigr).
$$
\end{proposition}

\begin{proof} Let $A_{K}^{sy}(2m,n):=\End_{\mbb_n}\bigl(V^{\otimes n}\bigr)$ be the
symplectic Schur algebra. By Theorem \ref{thm11}, there is a natural
surjection $KSp_{2m}(K)\twoheadrightarrow A_{K}^{sy}(2m,n)$. The
action of $KSp_{2m}(K)$ on $V^{\otimes n}/V^{\otimes n}\mbb_n^{(f)}$
factors through an action of $A_{K}^{sy}(2m,n)$. It follows that $$
\End_{KSp_{2m}(K)}\bigl(V^{\otimes n}/V^{\otimes
n}\mbb_n^{(f)}\bigr)=\End_{A_{K}^{sy}(2m,n)}\bigl(V^{\otimes
n}/V^{\otimes n}\mbb_n^{(f)}\bigr).
$$
Let $\overline{K}$ be the algebraic closure of $K$. By \cite{DDH}
and \cite{Oe},
$A_{K}^{sy}(2m,n)\otimes_{K}\overline{K}\cong
A_{\overline{K}}^{sy}(2m,n)$.
It follows that $$\begin{aligned}
&\quad\,\End_{KSp_{2m}(\overline{K})}\bigl(V_{\overline{K}}^{\otimes
n}/V_{\overline{K}}^{\otimes
n}\mbb_n^{(f)}\bigr)=\End_{A_{\overline{K}}^{sy}(2m,n)}\bigl(V_{\overline{K}}^{\otimes
n}/V_{\overline{K}}^{\otimes
n}\mbb_n^{(f)}\bigr)\\
&\cong\End_{A_{K}^{sy}(2m,n)}\bigl(V^{\otimes n}/V^{\otimes
n}\mbb_n^{(f)}\bigr)\otimes_{K}\overline{K}.
\end{aligned}$$
In particular, $$\begin{aligned} &\quad\,\dim_{K}\End_{KSp_{2m}(K)}\bigl(V^{\otimes n}/V^{\otimes
n}\mbb_n^{(f)}\bigr)=\dim_{K} \End_{A_{K}^{sy}(2m,n)}\bigl(V^{\otimes
n}/V^{\otimes
n}\mbb_n^{(f)}\bigr)\\
&=\dim_{\overline{K}}
\End_{\overline{K}Sp_{2m}(\overline{K})}\bigl(V_{\overline{K}}^{\otimes
n}/V_{\overline{K}}^{\otimes
n}\mbb_n^{(f)}\bigr)=\dim_{\overline{K}}\varphi''_{\overline{K}}(\mbb_n^{\overline{K}}).\end{aligned}
$$
Hence $\dim_{K}\End_{KSp_{2m}(K)}\bigl(V^{\otimes n}/V^{\otimes
n}\mbb_n^{(f)}\bigr)$ is independent of the infinite field $K$. Note
that
$\varphi''_{\overline{K}}(\mbb_n^{\overline{K}})=\varphi_{K}(\mbb_n)\otimes_{K}\overline{K}$.
It follows that $$
\dim_K\varphi_{K}(\mbb_n)=\dim_{\overline{K}}\varphi''_{\overline{K}}(\mbb_n^{\overline{K}})=\dim_{K}\End_{
KSp_{2m}(K)}\bigl(V^{\otimes n}/V^{\otimes n}\mbb_n^{(f)}\bigr),
$$
from which the proposition follows immediately.
\end{proof}

If $f=1$ and $m\geq n$, then it is easy to check that
$$\begin{aligned} \dim K\bBS_n&=\dim\mmC\bBS_n=\sum_{\lam\vdash n}(\dim
S_{\mmC}^{\lam})^2= \dim\End_{KSp_{2m}(\mmC)}\bigl(V_{\mmC}^{\otimes
n}/V_{C}^{\otimes n}\mbb_n^{(1)}\bigr)\\
&=\dim \End_{KG}\bigl(V^{\otimes n}/V^{\otimes n}\mbb_n^{(1)}\bigr),
\end{aligned}$$
This implies that $\varphi_K$ maps $K\bBS_n\cong \mbb_n/\mbb_n^{(1)}$
isomorphically onto $$\End_{KSp_{2m}(K)}\bigl(V^{\otimes
n}/V^{\otimes n}\mbb_n^{(1)}\bigr).$$ This gives the following corollary.

\begin{corollary}\label{maincor2} {\rm (\cite{DS}, \cite{Ma2})} Let $K$ be an arbitrary infinite field. If $m\geq n$, then
$\varphi_{1,K}$ maps $\mbb_n/\mbb_n^{(1)}\cong K\bBS_n$ isomorphically
onto $\End_{KSp(V)}\Bigl(V^{\otimes n}/V^{\otimes
n}\mbb_n^{(1)}\Bigr)$.
\end{corollary}

\begin{remark} Let $f$ be an integer with $1\leq f\leq
[n/2]$ and $K$ an infinite field. Proposition \ref{mainprop} proves
one side of the Brauer-Schur-Weyl duality between
$\mbb_n/\mbb_n^{(f)}$ and $KSp(V)$ on $V^{\otimes n}/V^{\otimes
n}\mbb_n^{(f)}$. We conjecture that the other side of the
Brauer-Schur-Weyl duality is also true. That is, the dimension of
the endomorphism algebra $\End_{\mbb_n/\mbb_n^{(f)}}\bigl(V^{\otimes
n}/V^{\otimes n}\mbb_n^{(f)}\bigr)$ is independent of $K$ and the
natural $K$-algebra homomorphism $\psi_{K}:
KSp_{2m}(K)\rightarrow\End_{\mbb_n/\mbb_n^{(f)}}\bigl(V^{\otimes
n}/V^{\otimes n}\mbb_n^{(f)}\bigr)$ is also surjective.
\end{remark}

\bigskip\bigskip

\end{document}